\documentclass[10pt]{article}
\usepackage[utf8]{inputenc}
\usepackage[T1]{fontenc}
\usepackage{amsfonts}
\usepackage{amsmath,latexsym}
\usepackage{amssymb,amsthm}
\usepackage{graphicx}
\usepackage[english]{babel}
\usepackage{hyperref}
\usepackage{color}
\usepackage{dsfont}
\usepackage{geometry}
\usepackage[french]{layout}
\theoremstyle{definition} 
\theoremstyle{remark}
\newtheorem{exam}{\textsf{Exemple}}[section]%[chapter]
\newtheorem{rmq}{\textsf{Remark}}%[section]%[chapter]

\theoremstyle{plain}
\newtheorem{theo}{\textsc{Theorem}}[section]%[chapter]
\newtheorem{Lemme}{\textsc{Lemma}}[section]%[chapter]
\newtheorem{prop}{\textsc{Proposition}}[section]%[chapter]
\newcommand{\cqfd}
{%
\mbox{}%
\nolinebreak%
\hfill%
\rule{2mm}{2mm}%
\medbreak%
\par%
}

\newcommand{\mc}[1]{\mathcal{#1}}

\usepackage{amsopn}

\DeclareMathOperator{\Supp}{Supp}
\DeclareMathOperator{\Div}{div}

\numberwithin{equation}{section}

\author{ 
	\textsc{Michel Duprez}\thanks{Inria, Universit\'e de Strasbourg,  ICUBE, \'equipe MIMESIS, Strasbourg, France. E-mail: {\tt
michel.duprez@inria.fr}.} \and \textsc{
	Pierre Lissy} \thanks{CEREMADE,  Universit\'e Paris-Dauphine \& CNRS UMR 7534, Universit\'e PSL, 75016 Paris, France. E-mail: {\tt
lissy@ceremade.dauphine.fr}.}}
	
\date{\today}

\title{Bilinear local controllability to the trajectories of the Fokker-Planck equation with a localized control}

\begin{document}
\maketitle

\begin{abstract}
This work is devoted to the control of the Fokker-Planck equation, posed on a smooth bounded domain of $\mathbb R^d$, with a localized drift force.
We prove that this equation is locally controllable to regular nonzero trajectories. Moreover, under some conditions, we explain how to reduce the number of controls around the reference control.  The results are obtained thanks to a standard linearization method and the fictitious control method. The main novelties are twofold. First, the algebraic solvability is performed  and used directly on the adjoint problem. We then prove a new Carleman inequality for the heat equation with a space-time varying first-order term: the right-hand side is the gradient of the solution localized on an open subset. We finally give an example of regular trajectory around which the Fokker-Planck equation is not controllable with a reduced number of controls, to highlight that our conditions are relevant.
\end{abstract}


%\tableofcontents

\section{Introduction and main results}
\subsection{Introduction}
\hspace*{4mm} Let $T>0$ and let $\Omega$ be a bounded domain in $\mathbb{R}^d$ ($d\in\mathbb N^*$),
regular enough (for example of class $\mathcal{C}^\infty$). 
Denote by $Q_T:=(0,T)\times\Omega$ and $\Sigma_T:=(0,T)\times\partial\Omega$. 
We consider the following system
\begin{equation}\label{syst primal}
\left\{\begin{array}{lll}
\partial_ty&=\Delta y+\Div(u y)&\mbox{in } Q_T,\\
y&=0&\mbox{on }\Sigma_T,\\
y(0,\cdot)&=y^0&\mbox{in }\Omega,
\end{array}\right.
\end{equation}
where $y^0\in L^2(\Omega)$ is the initial data and 
$u=(u_1,...,u_d)\in L^\infty((0,T)\times\Omega)^d$ is the control.

 It is well-known (see for instance \cite[Theorem and Proposition 3.1]{roberto}) 
that for  every initial data  $y^0\in L^2(\Omega)$ and every control $u\in L^\infty((0,T)\times\Omega)^d$, there exists a unique solution $y$ to System \eqref{syst primal} in the space 
$W(0,T)$, where  
\begin{equation*}
W(0,T):=L^2((0,T),H^1_0(\Omega ))
\cap H^1((0,T),H^{-1}(\Omega ))\hookrightarrow\mc{C}^0([0,T];L^2(\Omega)).
\end{equation*}

Equation \eqref{syst primal}, introduced in \cite{MR1512678}, is called the Fokker-Planck equation. %In the case where 
When the Fokker-Planck equation is posed on the whole space $\mathbb R^d$, it is strongly related to the stochastic differential equation (SDE)

\begin{equation}\label{SDE}
\left\{\begin{array}{lll}
dX_t&=\sum_{i=1}^du_i(X_t)dt+dW_t&\mbox{in } (0,T)\times \mathbb R^d,\\
X(0,\cdot)&=X^0&\mbox{in }\mathbb R^d,
\end{array}\right.
\end{equation}
where $W_t$ is the standard multi-dimensional Brownian motion starting from $0$. System \eqref{SDE} describes the movement of a particule of negligible mass, with constant and isotropic diffusion, under the action of a force field $u=(u_1,\ldots, u_d)$.

Under some regularity conditions on the drift term $U$, it is well-known that, by the It\^o Lemma, the probability density function $p$ associated to \eqref{SDE} verifies 

\begin{equation}\label{pws}
\left\{\begin{array}{lll}
\partial_t p&=\frac{1}{2}\Delta p+\Div(u p)&\mbox{in } (0,T)\times\mathbb R^d,\\
p(0,\cdot)&=p^0&\mbox{ in }\mathbb{R}^d,
\end{array}\right.
\end{equation}
where $p^0$ is some initial probability density function (see \textit{e.g.} \cite[Section 5.3]{MR2118834}).
By definition of a probability measure, we have $p^0\geqslant 0$ a.e. and $\int_{\mathbb R^d} p^0=1$. %It is then very easy to prove that these properties are preserved during time:
Then, we can easily prove the preservation of these properties during the time:
 any solution $p$ of System \eqref{pws} verifies also $p(t,\cdot)\geqslant 0$ a.e. and $\int_{\mathbb R^d} p(t,\cdot)=1$, for any $t\in[0,T]$ and hence remains a probability measure. We refer to \cite{MR987631} for more  explanations on the  Fokker-Planck equation, notably in the case of nonlinear drift terms or non-constant and anisotropic diffusion.

However, in the case where we impose Dirichlet boundary conditions as in \eqref{syst primal}, the derivation of the Fokker-Planck equation from a SDE is more %difficult:
complicated: the Brownian motion has to be replaced by an ``absorbed'' or ``killed'' Brownian motion, see \textit{e.g.} \cite[pp. 31-60]{MR1329992}. Moreover, the total mass of the initial condition is not conserved anymore, meaning that the probability of  remaining inside $\Omega$ decreases in time, and the solution to \eqref{syst primal} is not a probability density function anymore. We refer to \cite[Section 2]{roberto} for a discussion on the relevance of Dirichlet boundary conditions in this context. Neumann boundary conditions (that would restore the conservation of mass) \textcolor{black}{are beyond} the scope of the present article \textcolor{black}{(see the last item of Remark \ref{la4} for more explanations)}.

%While the 
 The controllability properties of the scalar linear heat equation in the case of a distributed control on an open subset and Dirichlet boundary condition are now well-understood (see notably \cite{lebeaurobbiano95} and \cite{fursikov1996controllability}). %,
The bilinear controllability seems to have been less explored. %Equation 
 The equation \eqref{syst primal} has been studied in \cite{MR1166481}, in the whole space and with controls localized everywhere in space and time. Concerning bilinear control %in the case where 
when the bilinear term $\Div(u y)$ is replaced by $uy$ with $u\in L^\infty((0,T)\times\Omega)$, we refer to \cite{MR3698164,MR2679642,MR2009955,MR1972539,MR1925029,MR2641453,MR2221515,MR3449653,MR3436584,MR2417139}.

%Let us mention that bilinear optimal control 
Optimal bilinear control of parabolic equations has previously been studied. 
 A first result was proved in \cite{MR1896168}, where a close %forth-order
 fourth-order in time model is investigated, with controls depending only on time. This result has been extended to second-order parabolic equations firstly in \cite{MR2724518} in the one-dimensional case,  then in \cite{MR2966923} in the multi-dimensional case, still for time-varying controls. For equation \eqref{syst primal} (in a slightly more general form), the case of space and time-varying  controls is treated in  \cite{roberto}. Notably, for a drift term that is affine in the control, the authors  prove the existence of optimal controls for general cost functionals, and derive first-order necessary optimality conditions using an adjoint state.
The controllability of the continuity equation, \textit{i.e.} System \eqref{syst primal} without diffusion, has been investigated in \cite{DMR19,DMR19min}.

%The  structure of the article is as follows. In
The paper is organized as follows: in
 Section \ref{s:m}, we %give 
 present the main results of the article (Theorem \ref{theo: contr traj 1}, resp. Theorem \ref{theo: contr traj 2}, which %gives
 provides a result of local controllability to the trajectories with $d$ \textcolor{black}{controls}, resp. a reduced number of controls around the reference control) and some remarks. Section \ref{s:li} is devoted to studying a linearized version of \eqref{syst primal}. In Section \ref{s:c}, we prove a new Carleman estimate (Proposition \ref{prop ine obs 1}) for solutions of the linear backward heat equation with \textcolor{black}{first-order terms}. The main novelty is that the local observation term is the gradient of the solution of the adjoint problem \eqref{sys:dual lin}%. This 
, which has already been proved in \cite{ML15} for constant coefficients. Moreover, %we are able to put 
we can put as many derivatives as we want in the left-hand side of our Carleman estimate, which will be %needed 
need for the rest of the proof. In Section \ref{s:a}, we explain how to  remove some components of the gradient  in the Carleman inequality. %This is performed by using 
\textcolor{black}{To demonstrate that, we use }we call an argument of  ``algebraic solvability''  (as introduced in \cite{92mcss} in the context of the stabilization of ODEs and in \cite{coronlissy2014} for the study of coupled systems of PDEs), based on ideas developed by Gromov in \cite[Section 2.3.8]{Gromovbook}. This procedure has already been used successfully in \cite{ACO,ML15,MR3820418,CG16,MR3624931,MR4011675,MR4011676}. The main novelty compared to the existing literature is that the algebraic solvability is performed directly on the dual problem. Moreover, %we are able to get rid 
\textcolor{black}{we can get rid }the high order derivatives of the right in order to obtain the final Carleman estimate \eqref{ine obs 2}.
In Section \ref{s:rc}, we use some arguments coming from optimal control theory in order to derive from our observability inequality the existence of regular enough controls, with a special form, in appropriate weighted spaces.  In Section \ref{s:nl}, we go back to the nonlinear problem by using  a standard strategy coming from \cite{MR3023058} together with some adapted inverse mapping Theorem.
% \textcolor{red}{
To finish, in Section \ref{sec:CE}, we give an example of \textcolor{black}{a}  trajectory around which the local controllability does not hold with a reduced number of controls.
%}

\subsection{\textcolor{black}{Main} results}
\label{s:m}

Let $(\overline{y},\overline{u})$ be a trajectory of \eqref{syst primal}, \textit{i.e.} verifying
\begin{equation}\label{syst linearise}
\left\{\begin{array}{lll}
\partial_t\overline{y}&=\Delta \overline{y}+\Div(\overline{u}\overline{y})&\mbox{in } Q_T,\\
\overline{y}&=0&\mbox{on } \Sigma_T,\\
\overline{y}(0,\cdot)&=\overline{y}^0\in L^2(\Omega)\setminus \{0\}&\mbox{in }\Omega.
\end{array}\right.
\end{equation}

\subsubsection{Controls with \textit{d} components}

%\textcolor{black}{
We first state a result of local controllability to the trajectories \textcolor{black}{for} System \eqref{syst primal} with a control containing $d$ components:
%}
\begin{theo}\label{theo: contr traj 1}
%Let $m\in\mathbb N^*$ (with possibly $m<d$). 
Let $\omega$ be any nonempty open subset of $\Omega$.
Assume that the trajectory $(\overline y,\overline{u})$ with $\overline{u}=(\overline{u}_1,...,\overline{u}_d)$  of System \eqref{syst linearise} is regular enough (for example of class $C^\infty$ on $(0,T)\times\Omega$).
%Assume that there exists  $q\in\mathbb N$ and some $(t_0,x_0)\in (0,T)\times\omega$ such that 
%\begin{equation}\label{cdet} \mbox{Rank}(M_q(\overline{u})(t_0,x_0))=d.
%\end{equation}
Then, System \eqref{syst primal} is locally controllable with localized controls, in the following sense: 

for every $\varepsilon>0$ and every $T>0$, there exists $\eta>0$ such that for any $y^0\in L^2(\Omega)$ verifying 
\begin{equation}\label{cieta}||y^0-\textcolor{black}{\overline{y}^0}||_{L^2(\Omega)}\leqslant \eta,\end{equation} there exists a trajectory $(y,u)$ to System \eqref{syst primal} such that 

\begin{equation*}
%\label{ll} 
\left\{\begin{array}{ll}y(T)&=\overline{y}(T),
\\u&=\overline{u}+v\mbox{ for some }v\in L^\infty((0,T)\times\Omega)^d,
\\ \rm{Supp}\,(v)&\subset (0,T)\times\omega,
\\ ||v||_{L^\infty((0,T)\times\Omega)^d}&\leqslant \varepsilon,
\\||y-\overline y||_{W(0,T)}&\leqslant \varepsilon.
\end{array}\right.
\end{equation*}
\end{theo}

\begin{rmq}
\begin{itemize}
\item The regularity assumptions on $(\overline y,\overline u)$ can be improved, notably it is enough that the reference trajectory is $C^r$ for some $r\in\mathbb N^*$ large enough, on an open subset of $(0,T)\times\textcolor{black}{\omega}$. 
%\item Remark that if $B=I_d$ (\textit{i.e.} we control every component of the gradient of $u$), condition \eqref{cdet} is automatically verified for $q=0$, whatever $\overline u$ is.
%\item Condition \eqref{cdet} notably implies that $q$ has to be chosen large enough such that 
%$N(m,q)\geqslant d-m$.
\item If $y^0=0$, the only solution to \eqref{syst primal} is $y\equiv 0$, whatever $u$ is, so that the only reachable state at time $T$ is $0$. As a consequence, $\eta>0$ has notably to be chosen small enough such that $y^0\not =0$. 
%\textcolor{black}{Notice that we also cannot consider $\overline y^0=0$ in \eqref{syst linearise}: in this case, we would have $\overline y \equiv 0$, so that notably $\overline y(T)=0$. By backward uniqueness (see \cite{MR0338517}), if the solution $y$ of \eqref{syst primal} verifies $y(T)=0$, then $y\equiv 0$ and notably $y^0=0$, which is excluded.} 
\item From the results given in \cite{Aronso},  as soon as $y^0\geqslant 0$, then any trajectory to System \eqref{syst primal} remains non-negative
(see also \cite{roberto}). This fact differs from the usual linear heat equation with internal control (see \cite{MR3669834}).
\item \textcolor{black}{We can also remark that we do not assume any relation between the control
domain $\omega$ and the support of $\bar u$. In particular, they can be disjoint.}
%\item \textcolor{black}{In Section \ref{sec:CE1}, we give an example of trajectory with the support of $\overline u$ not included in a set $[0,T]\times\omega_{\overline u}$ and for which the local controllability to the trajectories does not hold.}
%\item Assumption \eqref{cdet} is generic, in the following sense: if  $C^\infty((0,T)\times\omega)^2$ is endowed with the $C^q$ topology, the sets of the functions $(\overline y,\overline u\in C^\infty((0,T)\times\omega))^2$ verifying \eqref{cdet} is an dense open set.
 \end{itemize}
\end{rmq}

\subsubsection{Controllability acting through a control operator}

%\textcolor{black}{
In this section, we give a result of local controllability to the trajectories to System \eqref{syst linearise} with a control acting through a control operator $B\in\mathcal M_{d,m}(\mathbb R)$ with 
 $m\in\mathbb N^*$ such that $m\leqslant d$.
 %}
 
%\textcolor{black}{
We first introduce some notations.
%}
%Let $q\in \mathbb N$ and consider the following set 
%
%$$\mathcal E(m,q)= \{(\alpha_1,\ldots,\alpha_m)\in\mathbb N^m\ |\  0<\alpha_1+\ldots +\alpha_m\leqslant q \},$$
%with the convention that $\mathcal E(m,q)=\emptyset$ if $q=0$.
%Note that by an elementary computation, $$\#\mathcal E(m,q)=\frac{(q+m)!}{m!q!}-1=:N(m,q).$$
%%\textcolor{red}{pourquoi -1?}
%Let $m\in\mathbb N^*$ such that $m\leqslan:t d$ and let $B\in\mathcal M_{d,m}(\mathbb R)$. 
For $j\in\mathcal \{1,...,m\}$, we call $B_j^* \in \mathbb R^d$ the $j$-th line of $B^*$, and 
$$(B^*_{\textcolor{black}{j}}\cdot\nabla):\psi\in C^\infty(\mathbb R^d,\mathbb R)\mapsto B^*_j(\nabla \psi)\in C^\infty(\mathbb R^d,\mathbb R).$$
For $(\alpha_1,\ldots \alpha_m)\in \mathbb{N}^m$, %\mathcal E(m,q)$, 
we introduce the following operator:
\begin{equation*}\begin{aligned}
(B^*\cdot\nabla)^{\alpha_1,\ldots\alpha_m}:\psi\in C^\infty(\mathbb R^d,\mathbb R)\mapsto& \underbrace{(B^*_1\cdot\nabla)\ldots (B^*_1\cdot\nabla)}_{\alpha_1\mbox{ times}}\ldots\\&\underbrace{(B^*_m\cdot\nabla)\ldots(B^*_m\cdot\nabla)}_{\alpha_m\mbox{ times}}\psi\in C^\infty(\mathbb R^d,\mathbb R),\end{aligned}
\end{equation*}
%
%We introduce the following matrix:
%\begin{equation}\label{mtm}
%M_q(\overline{u})=\begin{pmatrix}&B^*_1&\\&\vdots& \\&B^*_m&\\ (B^*\cdot\nabla)^{1,0,\ldots,0}\overline u_1&\ldots &(B^*\cdot\nabla)^{1,0,\ldots,0}\overline u_d\\ 
%(B^*\cdot\nabla)^{0,1,\ldots,0}\overline u_1&\ldots &(B^*\cdot\nabla)^{0,1,\ldots,0}\overline u_d\\ 
% \vdots &\ddots&\vdots\\
%(B^*\cdot\nabla)^{0,\ldots,0,q}\overline u_1&\ldots  &(B^*\cdot\nabla)^{0,\ldots,0,q}\overline u_d\\ 
%\end{pmatrix}\in\mathcal M_{N(m,q)+m,d}(\mathbb R).
%\end{equation}
\textcolor{black}{and the family of $\mathbb R^d$ given by
\begin{equation*}
%\label{mtm}
M (\bar u)(t, x) = \{B_1^* , ..., B_m^* \} \cup \{((B^*\cdot\nabla)^\alpha \bar u_i(t, x))_{i\in \{1,\cdots,d\}}, \alpha\in\mathbb{N}^m,\alpha\neq 0\} .
\end{equation*}}

We have the following controllability result.
\begin{theo}\label{theo: contr traj 2}
Let $m\in\mathbb N^*$ (with possibly $m<d$). 
%Let $\omega$ be any nonempty open subset of $\Omega$.
%Assume that the trajectory $(\overline y,\overline{u})$ with $\overline{u}=(\overline{u}_1,...,\overline{u}_d)$  of System \eqref{syst linearise} is smooth enough (for example of class $C^\infty$ on $(0,T)\times\Omega$), and that there exists some open subset $\omega_{\overline u}$, strongly included in $\Omega$, such that  the support of $\overline u$ is included in $[0,T]\times\omega_{\overline u}$. 
Under the hypothesis of Theorem \ref{theo: contr traj 1}, assume that there exists  
%$q\in\mathbb N$ and 
some $(t_0,x_0)\in (0,T)\times\omega$ such that 
%\begin{equation}\label{cdet} \mbox{Rank}(M_q(\overline{u})(t_0,x_0))=d.
%\end{equation}
\textcolor{black}{\begin{equation}
\label{cdet} \textcolor{black}{\mbox{rank}}(M(\overline{u})(t_0,x_0))=d.
\end{equation}}
Then, System \eqref{syst primal} is locally controllable with localized controls, in the following sense: 

for every $\varepsilon>0$ and every $T>0$, there exists $\eta>0$ such that for any $y^0\in L^2(\Omega)$  verifying 
\begin{equation*}%\label{cieta2}
||y^0-\textcolor{black}{\overline{y}^0}||_{L^2(\Omega)}\leqslant \eta,
\end{equation*} there exists a trajectory $(y,u)$ to System \eqref{syst primal} such that 

\begin{equation*}
%\label{ll2} 
\left\{\begin{array}{ll}y(T)&=\overline{y}(T),
\\u&=\overline{u}+Bv\mbox{ for some }v\in L^\infty((0,T)\times\Omega)^m,
\\ \rm{Supp}\,(v)&\subset (0,T)\times\omega,
\\ ||v||_{L^\infty((0,T)\times\Omega)^m}&\leqslant \varepsilon,
\\||y-\overline y||_{W(0,T)}&\leqslant \varepsilon.
\end{array}\right.
\end{equation*}
\end{theo}
\begin{rmq}%\label{rmq2}
\begin{itemize}
%\item The regularity assumptions on $(\overline y,\overline u)$ can be improved, notably it is enough that the reference trajectory is $C^r$ for some $r\in\mathbb N^*$ large enough, on an open subset of $((0,T)\times\omega_0)$. 
\item Remark that if $B=I_d$ (\textit{i.e.} we control every component of the gradient of $u$), condition \eqref{cdet} is automatically verified for $q=0$, whatever $\overline u$ is. %\textcolor{black}{
Hence Theorem \ref{theo: contr traj 2} contains the result given in Theorem \ref{theo: contr traj 1}. Thus we will only give a proof of Theorem \ref{theo: contr traj 2}.
%}
%\item Condition \eqref{cdet} notably implies that $q$ has to be chosen large enough such that 
%$N(m,q)\geqslant d-m$.
%\item If $y^0=0$, the only solution to \eqref{syst primal} is $y\equiv 0$, whatever $u$ is, so that the only reachable state at time $T$ is $0$. As a consequence, $\eta>0$ has notably to be chosen small enough such that $y^0\not =0$. 
%\item From the results given in \cite{Aronso},  as soon as $y^0\geqslant 0$, then any trajectory to System \eqref{syst primal} remains non-negative
%(see also \cite{roberto}). This fact differs from the usual linear heat equation with internal control (see \cite{MR3669834}).

%\item Assumption \eqref{cdet} is generic, in the following sense: if  $C^\infty((0,T)\times\omega)^2$ is endowed with the $C^q$ topology, the sets of the functions $(\overline y,\overline u)\in C^\infty((0,T)\times\omega)^2$ verifying \eqref{cdet} is an dense open set. 
\item %\textcolor{black}{
In Section \ref{sec:CE}, we give an example of trajectory which does not satisfy condition \eqref{cdet}  and for which the local controllability to the trajectories does not hold.
It highlights that Condition \eqref{cdet} is \textcolor{black}{relevant}.
Even if the authors think that Condition \ref{cdet} is \textcolor{black}{not} optimal, 
find\textcolor{black}{ing} a necessary and sufficient condition remains on open problem.
\item \textcolor{black}{We can also remark that condition \eqref{cdet} is local on $\omega$. Notably, contrary to Theorem \ref{theo: contr traj 1}, if $m<d$, we necessarily have that the control
domain $\omega$ and the support of $\bar u$ intersect.}
%}
\end{itemize}
\end{rmq}

\begin{exam}
We give an explicit example, %in order to explain 
\textcolor{black}{to explain} better condition \eqref{cdet}. Let us assume that we want to control only the $m(<n)$ first components of the gradient, \textit{i.e.} 
$$B=\begin{pmatrix}1&0&\ldots &0\\
0&1&\ddots&\vdots\\
 \vdots  &\ddots&\ddots &0 \\
  0&\ldots &0&1\\
  0&\ldots&\ldots &0 \\ 
   \vdots  && &\vdots\\
    0&\ldots&\ldots &0 
    \end{pmatrix}\in \mathcal M_{n,m}(\mathbb R).$$
\textcolor{black}{
\textcolor{black}{For $\alpha=(\alpha_1,\ldots \alpha_m)\in \mathbb{N}^m$, we have 
\begin{equation*}
(B^*\cdot\nabla)^{\alpha_1,\ldots\alpha_m}(\psi)=\partial_{x_1}^{\alpha_1}\ldots \partial_{x_m}^{\alpha_m}(\psi).
\end{equation*}}
We deduce that 
\begin{equation*}
M(\overline{u})(t,x)=\{\textcolor{black}{e_1,\ldots,e_m}\}\cup \{(\partial_{x_1}^{\alpha_1}\ldots \partial_{x_m}^{\alpha_m}\bar u_i(t, x))_{i\in \{1,\cdots,d\}}, \alpha\in\mathbb{N}^m,\alpha\neq 0\},
\end{equation*}
where the vector $e_i$  is the $i$-th element of the canonical basis of $\mathbb{R}^d$.
%
%
%\begin{equation*}
%M_q(\overline{u})=\begin{pmatrix}
%1&0&\ldots&\ldots&\ldots &\ldots&0\\
%0&1&0&&&&0\\ 
%\vdots  &\ddots&\ddots &\ddots && &\vdots\\ 
%0&\ldots &0&1&0&\ldots &0\\
% \partial_{x_{1}}\overline u_{1}& \partial_{x_{1}}\overline u_{2}&\ldots &\ldots&\partial_{x_{1}}\overline u_{m+1}&\ldots&\partial_{x_{1}}\overline u_d\\ 
%\vdots &&&&\vdots&&\vdots\\
% \partial_{x_{m}}\overline u_{1}& \ldots&\ldots &\ldots&\partial_{x_{m}}\overline u_{m+1}&\ldots&\partial_{x_{m}}\overline u_d\\ 
%\partial^2_{x_{1}^2}\overline u_{1}&\ldots&\ldots &\ldots&\partial^2_{x_{1}^2}\overline u_{m+1}&\ldots&\partial^2_{x_{1}^2}\overline u_d\\\
%\vdots &&&&\vdots&&\vdots\\
%\ \partial^{p}_{x_m^q}\overline u_{1}&\ldots&\ldots&\ldots &\partial^{q}_{x_m^q} \overline u_{m+1}&\ldots&\partial^{q}_{x_m^q} \overline u_{m+1}\end{pmatrix}\in\mathcal M_{N(m,q)+m,d}(\mathbb R).
%\end{equation*}
We observe that there exists $(t_0,x_0)\in (0,T)\times\omega$ such that  the rank of the family  $M(\overline{u})(t_0,x_0)$ is equal to  $d$ if and only if there exists $(t_0,x_0)\in (0,T)\times\omega$ such that the rank of the family 
\textcolor{black}{$$ \{(\partial_{x_1}^{\alpha_1}\ldots \partial_{x_m}^{\alpha_m}\bar u_i(t, x))_{i\in \{m+1,\cdots,d\}}, \alpha\in\mathbb{N}^m,\alpha\neq 0\}$$}
is  equal to $d-m$.}

%We observe that $M_q(\overline{u})$ is of maximal rank  $d$ if and only if the following matrix:
%
%\begin{equation*}
%\widetilde M_q(\overline{u})=\begin{pmatrix}
%\partial_{x_{1}}\overline u_{m+1}&\ldots&\partial_{x_{1}}\overline u_d\\ 
%\vdots&&\vdots\\
%\partial_{x_{m}}\overline u_{m+1}&\ldots&\partial_{x_{m}}\overline u_d\\ 
%\partial^2_{x_{1}^2}\overline u_{m+1}&\ldots&\partial^2_{x_{1}^2}\overline u_d\\ \vdots&&\vdots\\\partial^{q}_{x_m^q} \overline u_{m+1}&\ldots&\partial^{q}_{x_m^q} \overline u_{m+1}\end{pmatrix}\in\mathcal M_{N(m,q),d-m}(\mathbb R),
%\end{equation*}
%is of maximal rank $d-m$.
\end{exam}

\section{Null controllability of the linearized system}
\label{s:li}
In what follows, we always assume that the trajectory $(\overline y,\overline u)$ of \eqref{syst linearise} verifies  the hypothesis of Theorem \ref{theo: contr traj 1}.
Consider the following linear parabolic system
\begin{equation}\label{syst lin primal}
\left\{\begin{array}{lll}
\partial_ty&=\Delta y+\Div(\overline u y)+\Div(\theta u)&\mbox{in } Q_T,\\
y&=0&\mbox{on } \Sigma_T,\\
y(0,\cdot)&=y^0&\mbox{in }\Omega,
\end{array}\right.
\end{equation}
where $y^0\in L^2(\Omega)$  and 
$\theta\in \mathcal{C}^\infty(\overline{\Omega})$ is such that 
  \begin{equation}\label{deftheta}
 \left\{\begin{array}{lll}
  \Supp(\theta)&\subseteq \omega,&\\
  \theta&\equiv1&\mathrm{~in~} \omega_0,\\
  0&\leqslant \theta\leqslant 1 &\mathrm{~in~} \Omega,
 \end{array}
\right.
\end{equation}
for some non-empty open subset $\omega_0$ which is strongly included in $\omega$.
The goal of this section is to prove the null controllability of System \eqref{syst lin primal}, with less controls than equations and regular enough controls in a special form.

\begin{rmq}\label{eqt}
Notice that the null controllability of \eqref{syst lin primal} is equivalent to the null controllability of the ``real'' linearized version of \eqref{syst primal} around $(\overline y,\overline u)$ given by 

\begin{equation}
\label{reale}
\left\{\begin{array}{lll}
\partial_ty&=\Delta y+\Div(\overline u y)+\Div(\overline y \tilde u)&\mbox{in } Q_T,\\
y&=0&\mbox{on } \Sigma_T,\\
y(0,\cdot)&=y^0&\mbox{in }\Omega.
\end{array}\right.
\end{equation}
\textcolor{black}{
Indeed,  since the solution of $(\overline y,\overline u)$ of \eqref{syst linearise} is in $C^\infty((0,T)\times\Omega)$, as soon as $\overline{y}^0\not = 0$,  on $(0,T)\times\omega$, $\overline y^{-1}(\{0\})$ is a closed subset of $(0,T)\times\omega$, which cannot be $(0,T)\times\omega$ since it has a finite $d$-dimensional Hausdorff measure in $\mathbb R^{d+1}$ (see \cite{MR1290401}).  Hence, $(0,T)\times\omega\setminus \overline y^{-1}(\{0\})$ contains a nonzero open subset, there exists some subset $(T_1,T_2)\times \tilde\omega$ of $(0,T)\times \omega$ such that $|\overline y|\geqslant C>0$ on $(T_1,T_2)\times\tilde\omega$, that we can assume to be exactly $(0,T)\times \omega$ without loss of generality. Hence, for any $i\in \{1,\ldots,d\}$, one can solve (in $\tilde u_i$) the equation $\theta u_i=\overline y \tilde u_i$ by posing 
$$\tilde u_i=\frac{\theta u_i}{\overline y}.$$
Remark that $\tilde u_i$ enjoys the same regularity properties as $u_i$.}
\end{rmq}
\subsection{Carleman estimates}
\label{s:c}
 Let us consider the following adjoint system associated to System \eqref{syst lin primal}
\begin{equation}\label{sys:dual lin}
 \left\{\begin{array}{lll}
-\partial_t\psi&=\Delta \psi+\overline u\cdot\nabla\psi&\mbox{ in }~ Q_T,\\
 \psi&=0&\mbox{ on }~\Sigma_T,\\
   \psi(T,\cdot)&=\psi^0&\mbox{ in }~\Omega.
        \end{array}
\right.
\end{equation}

First of all, we will introduce some notations. We denote by $|\cdot |$ the euclidean norm on $\mathbb R^M$, whatever $M\in\mathbb N^*$ is.
 For $s,\lambda>0$ and $p\geqslant 1$, let us define the two following functions:
\begin{equation}\label{defax}
 \alpha(t,x):=\dfrac{\exp((2p+2)\lambda\|\eta^0\|_{\infty})
 -\exp[\lambda(2p\|\eta^0 \|_{\infty}+\eta^0(x) )]}{t^p(T-t)^p}
\end{equation}
and
\begin{equation}\label{defax2}
 \xi(t,x):=\dfrac{\exp[\lambda(2p\|\eta^0 \|_{\infty}+\eta^0(x) )]}{t^p(T-t)^p}.
\end{equation}
Here, $\eta^0 \in\mathcal{C}^{\infty}(\overline{\Omega})$ is a function satisfying
\begin{equation*}
 |\nabla\eta^0 |\geqslant \kappa\mathrm{~in~}\Omega\backslash\omega_{1},~~~
  \eta^0>0\mathrm{~in~}\Omega ~~~\mathrm{and}~~~
   \eta^0=0\mathrm{~on~}\partial\Omega,
\end{equation*}
with $\kappa>0$ and $\omega_{1}$ some open subset verifying $\omega_{1}\subset\subset\omega_0$. The proof of the existence of such a function $\eta^0 $  
can be found  in \cite[Lemma 1.1, Chap. 1]{fursikov1996controllability} 
(see also \cite[Lemma 2.68, Chap. 2]{coron2009control}). 
We will use the two notations 
\begin{equation}\label{defaxs}
 \alpha^*(t):=\max\limits_{x\in\overline{\Omega}}\alpha(t,x)
 \mathrm{~~~and~~~}\textcolor{black}{\xi_*(t)}:=\min\limits_{x\in\overline{\Omega}}\xi(t,x),
\end{equation}
for all $t\in (0,T)$. Notice that these maximum and minimum are reached at the boundary $\partial\Omega$.
 For $s,\lambda>0$, let us define
\begin{equation}\label{defI}
 I(s,\lambda;u):=
 s^3\lambda^4\displaystyle\iint_{Q_T} e^{-2s\alpha}\xi^3u^2
 + s\lambda^2\displaystyle\iint_{Q_T} e^{-2s\alpha}\xi|\nabla u|^2.
\end{equation}

Let us now give some useful auxiliary results that we will need in our proofs. 
The first one is  a Carleman estimate which holds for solutions of the heat equation with 
non-homogeneous Neumann boundary conditions:

\begin{Lemme}\label{th carl neum} 

There exists a constant $C>0$ such that for any $u^0\in L^2(\Omega)$, $f_1\in L^2(Q_T)$ and $f_2\in L^2(\Sigma_T)$,
  the solution to the system 
 \begin{equation*}\left\{\begin{array}{lll}
  -\partial_tu-\Delta u &=f_1&\mathrm{in}~Q_T,\\
  \frac{\partial u}{\partial n}& =f_2&\mathrm{on}~\Sigma_T,\\
  u(T,\cdot)&=u^0&\mathrm{in~}\Omega\end{array}\right.
 \end{equation*}
satisfies
\begin{equation*}\begin{array}{r}
 I(s,\lambda;u)\leqslant C\left(
 s^3\lambda^4\displaystyle\iint_{(0,T)\times\omega_{1}} e^{-2s\alpha}\xi^3u^2
 +s\lambda\displaystyle\iint_{\Sigma_T}e^{-2s\alpha^*}\xi_*f_2^2 
 ~~~~~~~~~~~~~~~~~~~~~\right.\\\left.
 +\displaystyle\iint_{Q_T} e^{-2s\alpha}f_1^2\right),
\end{array}\end{equation*}
for all $\lambda\geqslant C$ and $s\geqslant C(T^p+T^{2p})$.
\end{Lemme}
Lemma \ref{th carl neum} is proved in \cite[Theorem 1]{GuerreroFourier} in the case $p=1$. However, following the steps of the proof given in \cite{GuerreroFourier}, one can prove exactly the same inequality for any $p\in\mathbb N^*$.

From Lemma \ref{th carl neum}, one can deduce the following result:

\begin{Lemme}\label{carleman 2X2 constant}
 Let $f\in L^2(\Sigma_T)$, $G=(g_1,\ldots g_d)\in L^\infty(Q_T)^d$ and $h\in L^2(Q_T)$. % and $a\in L^\infty(Q_T)$.
Then, there exists a constant $C>0$ such that 
for every $\varphi^T\in L^2(\Omega)$, 
the solution $\varphi$ 
 to the system 
\begin{equation*}
%\label{primal 2x2 constant fourier}
 \left\{\begin{array}{lll}
-\partial_t\varphi&=\Delta \varphi+G\cdot\nabla\varphi%+a\varphi
+h&\mathrm{in}~ Q_T,\\
       \frac{\partial\varphi}{\partial n}&=f&\mathrm{on}~\Sigma_T,\\
       \varphi(T,\cdot)&=\varphi^T&\mathrm{in}~\Omega
        \end{array}
\right.
\end{equation*}
satisfies
\begin{equation*}\begin{array}{ll}
 I(s,\lambda;\varphi)&
 \leqslant C\left(
 s^3\lambda^4\displaystyle\iint_{(0,T)\times\omega_{1}} e^{-2s\alpha}\xi^3\varphi^2
 %\right.\\
 %~~~~~~~~~~~~~~~~~~~~~~~~~~~~~~~~~~~~~~~~~~~~~~~~~~~~~~~~~~~~~~~~~
 %\left.
 +s\lambda\displaystyle\iint_{\Sigma_T}e^{-2s\alpha^*}\xi_*f^2 \right .
  \\&\left .+\displaystyle\iint_{Q_T}e^{-2s\alpha}h^2 \right),
\end{array}\end{equation*}
for every $\lambda\geqslant C$ and $s\geqslant s_0=C(T^p+T^{2p})$.
\end{Lemme}
The proof of Lemma \ref{carleman 2X2 constant} is standard and is left to the reader (one just has to apply Lemma \ref{th carl neum} 
and absorb the remaining lower-order terms thanks to the left-hand side).

We will also need the following estimates.

\begin{Lemme}\label{poincare poids}
Let $r\in \mathbb{R}$. Then, there exists $C:=C(r,\omega_1,\Omega)>0$  
such that, for every $T>0$ and every $u\in L^2((0,T),H^1(\Omega))$,
\begin{equation*}\begin{array}{r}
 s^{r+2}\lambda^{r+2}\displaystyle\iint_{Q_T} e^{-2s\alpha}\xi^{r+2}u^2
 \leqslant C\left(s^r\lambda^r
  \displaystyle\iint_{Q_T} e^{-2s\alpha}\xi^{r}|\nabla u|^2
  \right.~~~~~~~~~~~~~~~~~~~~~~~\\\left.
  + s^{r+2}\lambda^{r+2}\displaystyle\iint_{(0,T)\times\omega_1} e^{-2s\alpha}\xi^{r+2}u^2\right),
\end{array}\end{equation*}
for every $\lambda\geqslant C$ and $s\geqslant C(T^{2p})$.
\end{Lemme}

The proof of this lemma can be found for example in \cite[Lemma 3]{coronguerrero2009} in the case $p=9$. However, following the steps of the proof given in \cite{coronguerrero2009}, one can prove exactly the same inequality for any $p\in\mathbb N^*$.

%In order to deal
\textcolor{black}{To deal} with more regular solutions, one needs the following lemma.

\begin{Lemme} \label{lemme regul}
Let  $z_0\in H^1_0(\Omega)$, $G\in C^{\infty}(Q_T)^d$ and $f\in  L^2(Q_T)^m$. 
Let us denote by $\mathcal{R}:=-\Delta-G\cdot\nabla$ and 
  consider the solution $z$ to the system
    \begin{equation*}
    %\label{systeme regul}
 \left\{\begin{array}{lll}
\partial_tz&=\Delta z+G\cdot\nabla z
+f&\mathrm{in}~ Q_T,\\
       z&=0&\mathrm{on}~\Sigma_T,\\
       z(0,\cdot)&=z_0&\mathrm{in}~\Omega.
        \end{array}
\right.
\end{equation*}
Let $n\in\mathbb N$. Let us assume that 
$z_0\in H^{2n+1}(\Omega)$, $f\in L^2((0,T),H^{2n}(\Omega))\cap H^{n}((0,T),L^{2}(\Omega))$ 
and satisfy the following compatibility conditions:
\begin{equation}\label{coco}
 \left\{\begin{array}{l}
      g_0:=z_0\in H^1_0(\Omega),  \\
      g_1:= f(0,\cdot)-\mathcal{R}g_0\in H^1_0(\Omega), \\
      \vdots\\
      \textcolor{black}{g_{n}}:= \partial^{n-1}_t f(0,\cdot)-\mathcal{R}\textcolor{black}{g_{n-1}}\in H^1_0(\Omega).
        \end{array}
\right.
\end{equation}
 Then $z\in L^2((0,T),H^{2n+2}(\Omega))\cap H^{n+1}((0,T),L^{2}(\Omega))$
  and we have the estimate
 \begin{equation*}\begin{aligned}
 %\label{estim regul 1}
&  \|z\|_{ L^2((0,T),H^{2n+2}(\Omega))\cap H^{n+1}((0,T),L^{2}(\Omega))}
 \\& \leqslant C(\|f\|_{L^2((0,T),H^{2n}(\Omega))\cap H^n((0,T),L^{2}(\Omega))}
  +\|z_0\|_{H^{2n+1}(\Omega)}).\end{aligned}
 \end{equation*}

\end{Lemme}
It is a classical result that can be easily deduced for example from \cite[Th. 6, p. 365]{MR2597943}.

We are now able to prove the following crucial inequality:
\begin{prop}\label{prop ine obs 1}
 Let \textcolor{black}{$\mu>0$ and} $N\in\mathbb N$ with $N\geqslant 3$
 % and assume that for any $i\in \{1,...,d\}$, $$\overline u_i\in W^{\lfloor \frac{N-1}{2}\rfloor,\infty}((0,T),W^{N,\infty}(\Omega)).$$ 
 . Then, there exists \textcolor{black}{$p\geqslant 2$ and} $C>0$  such that for every 
  $\psi^0\in L^2(\Omega)$, 
 the corresponding solution $\psi$ to System \eqref{sys:dual lin} satisfies
\begin{equation}\label{ine obs 1}
\begin{array}{l}
\lambda^{2}\displaystyle\iint_{Q_T} e^{-2s\alpha-2\mu s\alpha^*}(s\xi)|\nabla^{N+1}\psi|^2 \\+\ldots+\lambda^{2N+2}\displaystyle\iint_{Q_T} e^{-2s\alpha-2\mu s\alpha^*}(s\xi)^{2N+1}|\nabla\psi|^2\\ + \lambda^{2N+2}\displaystyle\iint_{Q_T}e^{-2s\alpha^*-2\mu s\alpha^*}(s\xi_*)^{2N+1}|\psi|^2\\
\leqslant C
\lambda^{2N+2}\displaystyle\iint_{(0,T)\times\omega_0}  e^{-2s\alpha-2\mu s\alpha^*}\left(s\xi\right)^{2N+1}|\nabla\psi|^2
\end{array}
\end{equation}
for every $\lambda\geqslant C$ and $s\geqslant s_0=C(T^p+T^{2p})$.
\end{prop}

Such a Carleman inequality seems new to the authors in the context of non-constant coefficients (%it was proved
\textcolor{black}{proved} in \cite{ML15} in the case of constant coefficients). The main improvement comes from the fact that the observation is a gradient of the solution $\psi$ on $\omega_0$ (and not the solution itself). We are also able to introduce as many derivatives of $\psi$ as we want in the left-hand side, as soon as $\overline u_i$ is regular enough.

\begin{rmq}\label{la4}
\begin{itemize}
\item Notice that the proof proposed here relies on the fact that the lower-order terms in equation \eqref{sys:dual lin} are of \textcolor{black}{first order}, and would fail in the presence of lower-order terms of order $0$. Indeed, in the first step of our proof (inequality \eqref{ine: neum carl proof}), some term that cannot be absorbed will appear.
\item Notice that inequality \eqref{ine obs 1} automatically implies that any solution  $\psi$ of  \eqref{sys:dual lin} lives in high order weighted Sobolev spaces.
This is not a surprise since we know that away from the final time $t=T$, any solution of \eqref{sys:dual lin} is regular.
\item \textcolor{black}{Remark that the proof provided here would fail for Neumann boundary conditions, since the argument in our last step, based on a Poincar\'e-like inequality, is not true anymore. It is not clear for the authors how one can adapt it in this case.}
\end{itemize}
\end{rmq}

\textbf{Proof of Proposition \ref{prop ine obs 1}.}

The proof is inspired by \cite{coronguerrero2009} and is quite similar to \cite{ML15}. Let $\mu>0$. In all what follows, $C>0$ is a constant that does not depend on $s$ or $\lambda$ (but that might depend on the other parameters, notably $p$, $N$, $\eta$, $T$, $\mu$) and that might change from inequality to inequality. We assume without loss of generality that $N$ is odd (the case $N$ even can be treated similarly).

Let $\psi$ the solution to System  \eqref{sys:dual lin}.
We introduce the following auxiliary functions:
\begin{equation}
\label{rp1}
\rho_1^*:=e^{-\mu s\alpha^*},\mbox{    }\psi_1:=\rho_1^*\psi.\end{equation}
Then $\psi_1$ is solution of

\begin{equation}\label{sys:psi0}
 \left\{\begin{array}{lll}
-\partial_t\psi_1&=\Delta\psi_1
+\overline u\cdot \nabla\psi_1-\partial_t\rho_1^*\psi
&\mbox{ in }~ Q_T,\\
 \psi_1&=0&\mbox{ on }~\Sigma_T,\\
   \psi_1(T,\cdot)&=0&\mbox{ in }~\Omega.
        \end{array}
\right.
\end{equation}
We remark that  $\phi:=\nabla^{N}\psi_1$ (the operator $\nabla$ applied $N$ times, or in other words, all the derivatives of order $N$ of $\psi_1$, ordered for example lexicographically) satisfies the system 
\begin{equation*}%\label{primal 2x2 constant dirch}
 \left\{\begin{array}{ll}
-\partial_t\phi=\Delta\phi
+\sum\limits_{i=1}^{N}G_i\cdot \nabla^i\psi_1+\overline{u}\cdot\nabla\phi-\partial_t\rho_1^*\nabla^{N} \psi&\mbox{ in }~ Q_T,\\
 \frac{\partial\phi}{\partial n}= \frac{\partial\phi}{\partial n}&\mbox{ on }~\Sigma_T,\\
   \phi(T,\cdot)=0&\mbox{ in }~\Omega,
        \end{array}
\right.
\end{equation*}
where, for any $i\in \{1,...,N\}$, $G_i$ is   an  essentially bounded tensor 
of \textcolor{black}{appropriate} size, whose coefficients are depending only on
$\overline u_i$ and its derivatives in space up to the order $i$.
Applying Lemma \ref{carleman 2X2 constant} to the different components of $\phi$, we obtain the following estimate
\begin{equation}\label{ine: neum carl proof}
\begin{array}{l}
 I(s,\lambda;\phi)
 \leqslant C\left(
\underbrace{ s\lambda\displaystyle\iint_{\Sigma_T}e^{-2s\alpha^*}\xi_*
 \left| \frac{\partial \phi}{\partial n}\right|^2}_{(I)}+\underbrace{ \displaystyle\iint_{Q_T}e^{-2s\alpha}|\partial_t\rho_1^*\nabla^{N} \psi|^2 }_{(II)}
 \right.\\
 ~~~~~~~~~~~~~~~~~~~~~~~~~~~~~~~~~~~~~
 \left.

+\underbrace{ \displaystyle\iint_{Q_T}e^{-2s\alpha}\sum_{i=1}^{N}|\nabla^i \psi_1|^2 }_{(III)} +s^3\lambda^4\displaystyle\iint_{(0,T)\times\omega_1} e^{-2s\alpha}\xi^3|\phi|^2 \right).
\end{array}
\end{equation}
The rest of the proof is divided into four steps:
\begin{itemize}
\item[$\bullet$] In a first step, we will estimate the \textcolor{black}{ boundary term (I) }
%appearing in the right-hand side of \eqref{ine: neum carl proof} 
by some global interior term involving $\psi_1$,  which will be %absorb
 \textcolor{black}{absorbed} later on (in the last step). We will also absorb the \textcolor{black}{ term (II)}  under some condition on $p$.
\item[$\bullet$] In a second step, we will estimate the \textcolor{black}{ term (III)} by some local terms involving $\nabla \psi_1$ and its derivatives on $\omega_1$, and get rid of the third term of the right-hand side.
\item[$\bullet$] In a third step, we will estimate the high-order local terms created at the previous step by some local terms involving only $\nabla \psi_1$ on $\omega_0$.
\item[$\bullet$] In a last step, we will use some Poincar\'e-like inequality in order to recover the variable $\psi$ in the left-hand side and bound the global interior term of the right-hand side involving $\psi_1$ by an interior term involving $\nabla \psi$. We will conclude by coming back to the original variable $\psi$, in order to establish \eqref{ine obs 1}.

\end{itemize}

\noindent{\bf Step 1:}
Let  $\tilde\theta\in\mathcal{C}^2(\overline{\Omega})$ 
a function satisfying
\begin{equation*}
\frac{\partial\tilde\theta}{\partial n}=\tilde\theta=1  \mathrm{~on~} \partial\Omega.  
\end{equation*}
An integration by parts of the boundary term leads to
\begin{equation*}\begin{array}{l}
s\lambda\displaystyle\int_{0}^Te^{-2s\alpha^*}\xi_*\displaystyle\int_{\partial\Omega}
 \left|\frac{\partial \phi}{\partial n}\right|^2 
=s\lambda\displaystyle\int_{0}^Te^{-2s\alpha^*}\xi_*\displaystyle\int_{\partial\Omega}
 \frac{\partial\phi}{\partial n}\nabla\phi\cdot\nabla\tilde\theta  \\
=s\lambda\displaystyle\int_{0}^Te^{-2s\alpha^*}\xi_*\displaystyle\int_{\Omega}
 \Delta\phi\nabla\phi\cdot\nabla\tilde\theta  
+s\lambda\displaystyle\int_{0}^Te^{-2s\alpha^*}\xi_*\displaystyle\int_{\Omega}
 \nabla(\nabla\tilde\theta\cdot\nabla\phi)\cdot\nabla\phi .
 \end{array}\end{equation*}
Hence 
  \begin{equation*}%\label{estim preuve lemme carl3}
 \begin{array}{rcl}
s\lambda\displaystyle\int_{0}^Te^{-2s\alpha^*}\xi_*\displaystyle\int_{\partial\Omega}
 \left|\frac{\partial\phi}{\partial n}\right|^2 
\leqslant C\lambda\displaystyle\int_{0}^Te^{-2s\alpha^*}
 s\xi_*\|\psi_1\|_{H^{N+2}(\Omega)}\|\psi_1\|_{H^{N+1}(\Omega)}.
 \end{array}
 \end{equation*}
 Using the interpolation inequality
\begin{equation*}
\|\psi_1\|_{H^{N+2}(\Omega)}\leqslant C\|\psi_1\|_{H^{N+1}(\Omega)}^{1/2}\|\psi_1\|_{H^{N+3}(\Omega)}^{1/2}
\end{equation*} 
 and Young's inequality $ab\leqslant \frac{a^q}{q}+\frac{b^{q'}}{{q'}}$ ($\frac{1}{q}+\frac{1}{q'}=1$) for $a,b\geqslant 0$ and $q=4$, we deduce that for any $c\in \mathbb R$, we have 
\begin{equation}\label{estim preuve lemme carl3}
 \begin{array}{lll}
&\lambda\displaystyle\int_{0}^Te^{-2s\alpha^*}s\xi_*\displaystyle\int_{\partial\Omega}
 \left|\frac{\partial\phi}{\partial n}\right|^2 
\\&\leqslant C\lambda\displaystyle\int_{0}^Te^{-2s\alpha^*}
(s\xi_*)^{c} \|\psi_1\|_{H^{N+3}(\Omega)}^{1/2}(s\xi_*)^{(1-c)}\|\psi_1\|_{H^{N+1}(\Omega)}^{3/2} \\
 &\leqslant \textcolor{black}{C\lambda\displaystyle\int_{0}^Te^{-2(1+\mu)s\alpha^*}(s\xi_*)^{4c} \|\psi\|_{H^{N+3}(\Omega)}^2 }\\
%&&~~~~~~~~~~~~~~~~~~
 &\textcolor{black}{+C\lambda\displaystyle\int_{0}^Te^{-2(1+\mu)s\alpha^*}(s\xi_*)^{\frac{4(1-c)}{3}}
 \|\psi\|_{H^{N+1}(\Omega)}^2.}
 \end{array}\end{equation}
Consider the function $\textcolor{black}{\psi_2:=\rho^*_2\psi}$, where 
\textcolor{black}{\begin{equation}\label{rho1}\rho^*_{2}:= (s\xi_*)^{\frac{2(1-c)}{3}}e^{-(1+\mu)s\alpha^*}.\end{equation}}
The function $\psi_2$ is solution to the system
   \begin{equation*}
%   \label{dual 2x2 psi rond}
 \left\{\begin{array}{rll}
-\partial_t\psi_2&=\Delta \psi_2
+\overline u\cdot\nabla\psi_2
\textcolor{black}{-\partial_t(\rho_{2}^*)\psi}%-\rho_{1}^*\partial_t(\rho_0^*)\psi
&\mbox{ in }~ Q_T,\\
 \psi_2&=0&\mbox{ on }~\Sigma_T,\\
   \psi_2(T,\cdot)&=0&\mbox{ in }~\Omega.
        \end{array}
\right.
\end{equation*}
Using Lemma \ref{lemme regul}  for $\psi_2$ (remark that the compatibility conditions \eqref{coco} are verified,  since \textcolor{black}{$\psi_2(T,\cdot)=0$  and $\partial^j_t \rho_{2}^*(T,\cdot)=0$ for any $j\in \mathbb N$}), we deduce that  %\begin{array}{rcl}

 \begin{equation}\label{boot2}
 \begin{array}{ll}
&  \|\psi_2\|_{L^2((0,T),H^{2n+2}(\Omega))\cap H^{n+1}((0,T),L^2(\Omega))}\\
 & \leqslant C\|\textcolor{black}{\partial_t(\rho_{2}^*)\psi}\|_{ L^2((0,T),H^{2n}(\Omega))\cap H^{n}((0,T),L^2(\Omega))},
 \end{array}
 \end{equation}
 for $n=1,2,\ldots ,(N+1)/2$.
The definitions of $\xi_*$ and $\alpha^*$ given in \eqref{defaxs}, the definition of $\rho_2^*$ given in \eqref{rho1} lead to 
\begin{equation}\label{rhot}
 \textcolor{black}{
 %\left\{\begin{array}{ll}
 % |\partial_t^k\rho_{0}^*|&\leqslant C(s\xi_*)^{k+\frac{k}{p}}e^{-\mu s\alpha^*},\\
  |\partial_t^k\rho_{2}^*|\leqslant C(s\xi_*)^{\frac{2(1-c)}{3}+k+\frac{k}{p}}e^{-s(1+\mu)\alpha^*}
 %|\partial_t^k\rho_{1}^*|&\leqslant C(s\xi_*)^{\frac{2(1-c)}{3}+k+\frac{k}{p}}e^{-s\alpha^*}.
 % \end{array}\right.}
 }
 \end{equation}
 for $k\in\{1,\dots,\frac{N+3}{2}\}$ (we recall that $C$ can depend on $\mu$). 
%\begin{equation}\label{rhot}
% \left\{\begin{array}{ll}
%  |\partial_t\rho_{0}^*|&\leqslant C(s\xi_*)^{1+\frac{1}{p}}e^{-\mu s\alpha^*},\\
%  \vdots\\
%  |\partial^{\frac{N+3}{2}}_{t}\rho_{0}^*|&\leqslant C(s\xi_*)^{\frac{N+3}{2}+\frac{N+3}{2p}}e^{-\mu s\alpha^*},\\
% |\partial_t\rho_{1}^*|&\leqslant C(s\xi_*)^{\frac{2(1-c)}{3}+1+\frac{1}{p}}e^{-s\alpha^*},\\
% \vdots\\
%  |\partial^{\frac{N+3}{2}}_{t}\rho_{1}^*|&\leqslant C(s\xi_*)^{\frac{2(1-c)}{3}+\frac{N+3}{2}+\frac{N+3}{2p}}e^{-s\alpha^*}. \end{array}\right.
% \end{equation}
 Remark that for any $k\leqslant l$, we have 
 \begin{equation}\label{rhot2}
 \textcolor{black}{  |\partial^{k}_{t}\rho_{2}^*|\leqslant C|\partial^{l}_{t}\rho_{2}^*|.}
 \end{equation}
 Combining  \eqref{boot2} for $n=(N-1)/2$, \eqref{rhot}, \eqref{rhot2} and the equations satisfied by $\psi$ and $\psi_1$, we obtain
 \textcolor{black}{ \begin{multline}\label{estim43}
\lambda\displaystyle\int_{0}^Te^{-2(1+\mu)s\alpha^*}(s\xi_*)^{\frac{4(1-c)}{3}}
 \|\psi\|_{H^{N+1}(\Omega)}^2\\
 \leqslant C \lambda\left(\displaystyle\int_{0}^Te^{-2(1+\mu)s\alpha^*}(s\xi_*)^{\frac{4(1-c)}{3}+{N+1}+\frac{N+1}{p}}
||\psi||^2_{L^2(\Omega)} 
 \right. \\ \left.
+\displaystyle\int_{0}^Te^{-2(1+\mu)s\alpha^*}(s\xi_*)^{\frac{4(1-c)}{3}+{2}+\frac{2}{p}}||\psi||^2_{H^{N-1}(\Omega)} \right). 
\end{multline}}
In the right-hand side of \eqref{estim43}, we would like to estimate the term 
\textcolor{black}{$$ \displaystyle\int_{0}^Te^{-2(1+\mu)s\alpha^*}(s\xi_*)^{\frac{4(1-c)}{3}+{2}+\frac{2}{p}}
||\psi||^2_{H^{N-1}(\Omega)}  .$$
}
This can be done using exactly the same processus by introducing some appropriate auxiliary weight that multiplies $\psi$ or $\psi_2$ as in \eqref{rho1}, using Lemma \ref{lemme regul} successively for  $n=(N-1)/2,\ldots,0$, \eqref{rhot} and \eqref{rhot2}. At the end, by gathering all the inequalities, we obtain
  \textcolor{black}{ \begin{equation}\label{estim44}\begin{aligned}
\lambda\displaystyle\int_{0}^Te^{-2(1+\mu)s\alpha^*}(s\xi_*)^{\frac{4(1-c)}{3}}
 \|\psi\|_{H^{N+1}(\Omega)}^2 
\\ \leqslant C \lambda\displaystyle\int_{0}^Te^{-2(1+\mu)s\alpha^*}(s\xi_*)^{\frac{4(1-c)}{3}+{N+1}+\frac{N+1}{p}}
||\psi||^2_{L^2(\Omega)}  . \end{aligned}
\end{equation}}
Applying the same technique also leads to 
\textcolor{black}{   \begin{equation}\label{estim45}\begin{aligned}
\lambda\displaystyle\int_{0}^Te^{-2(1+\mu)s\alpha^*}(s\xi_*)^{4c}
 \|\psi\|_{H^{N+3}(\Omega)}^2 \\
 \leqslant C \lambda\displaystyle\int_{0}^Te^{-2(1+\mu)s\alpha^*}(s\xi_*)^{4c+{N+3}+\frac{N+3}{p}}
||\psi||^2_{L^2(\Omega)}  
. \end{aligned}
\end{equation}}
 From \eqref{estim preuve lemme carl3}, \eqref{estim44} and \eqref{estim45}, we deduce that 
  \textcolor{black}{\begin{multline}\label{estim preuve lemme carl4}
\lambda\displaystyle\int_{0}^Te^{-2(1+\mu)s\alpha^*}s\xi_*\displaystyle\int_{\partial\Omega}
 \left|\frac{\partial\phi}{\partial n}\right|^2 
 \\\leqslant C \lambda\left(
\displaystyle\int_{0}^Te^{-2(1+\mu)s\alpha^*}(s\xi_*)^{\frac{4(1-c)}{3}+{N+1}+\frac{N+1}{p}}
||\psi||^2_{L^2(\Omega)}  \right .
\\\left . +\displaystyle\int_{0}^Te^{-2(1+\mu)s\alpha^*}(s\xi_*)^{4c+{N+3}+\frac{N+3}{p}}
||\psi||^2_{L^2(\Omega)}  
\right). \end{multline} }
 Since we would like the powers in the right-hand side to be equal, it is natural to impose that 
 $$ 4c+{N+3}+\frac{N+3}{p}=\frac{4(1-c)}{3}+{N+1}+\frac{N+1}{p},$$
\textit{i.e.} 
 \begin{equation}\label{condc} 
 c=\frac{-3-p}{8p}. 
 \end{equation}
Thus, using \eqref{estim preuve lemme carl4} and \eqref{condc},
 we deduce that  
\textcolor{black}{ \begin{equation}\label{estim bord2}
\lambda\displaystyle\int_{0}^Te^{-2s\alpha^*}s\xi_*\displaystyle\int_{\partial\Omega}
 \left|\frac{\partial\phi}{\partial n}\right|^2 
 \leqslant C \lambda\displaystyle\int_{0}^Te^{-2(1+\mu)s\alpha^*}(s\xi_*)^{\frac{2 N (p+1)+5 p+3}{2 p}}
||\psi||^2_{L^2(\Omega)}  .\end{equation}}
From \eqref{ine: neum carl proof}, \eqref{estim bord2}, the first line of \eqref{rhot} and the definition of $\psi_1$ given in \eqref{rp1}, we already deduce that 
\begin{equation*}
%\label{i1}
\begin{array}{l}
 I(s,\lambda;\phi)
 \leqslant C\left(
 s^3\lambda^4\displaystyle\iint_{(0,T)\times\omega_1} e^{-2s\alpha}\xi^3|\nabla^N \psi_1|^2 \right . \\ \left .+ \lambda 
 \displaystyle\iint_{Q_T}e^{-2s\alpha^*}(s\xi_*)^{\frac{2 N (p+1)+5 p+3}{2 p}}
|\psi_1|^2   \right. \\ \left . +\displaystyle\iint_{Q_T}e^{-2s\alpha}\sum_{i=1}^{N}|\nabla^i \psi_1|^2  +\displaystyle\iint_{Q_T}e^{-2s\alpha}(s\xi_*)^{2+\frac{2}{p}}|\nabla^N \psi_1|^2 \right).
\end{array}
\end{equation*}
By definition  of $\xi_*$ given in \eqref{defaxs}, it is clear that $\xi_*\leqslant \xi$. Hence, taking $p$ large enough such that $2+\frac{2}{p}\leqslant 3$ (\textit{i.e.} $p\geqslant 2$), $s,\lambda$ large enough and using the definition of $ I(s,\lambda;\phi)$ given in \eqref{defI}, we deduce that we can absorb the last term of the right-hand-side, so that we obtain
\begin{equation}\label{i1}
\begin{array}{l}
 I(s,\lambda;\phi)
 \leqslant C\left(
 s^3\lambda^4\displaystyle\iint_{(0,T)\times\omega_1} e^{-2s\alpha}\xi^3|\nabla^N \psi_1|^2 \right . \\ \left .+ \lambda 
 \displaystyle\iint_{Q_T}e^{-2s\alpha^*}(s\xi_*)^{\frac{2 N (p+1)+5 p+3}{2 p}}
|\psi_1|^2   \right. \\ \left . +\displaystyle\iint_{Q_T}e^{-2s\alpha}\sum_{i=1}^{N}|\nabla^i \psi_1|^2 \right).
\end{array}
\end{equation}
{\bf Step 2: }
We apply Lemma \ref{poincare poids} successively with 
$$(u,r)=(\nabla^{N-1} \psi_1,3), \ldots,(u,r)=(\nabla\psi_1,2N-1).$$ 
We obtain a sequence of inequalities of the form
\begin{equation*}
%\label{gauche1}
\begin{array}{l}
 s^{5}\lambda^{6}\displaystyle\iint_{Q_T} e^{-2s\alpha}\xi^{5}|\nabla^{N-1}\psi_1|^2  \leqslant C\left(s^{3}\lambda^{4}\displaystyle\iint_{Q_T} e^{-2s\alpha}\xi^{3}|\nabla^{N}\psi_1|^2  \right . \\ \left .+ s^{5}\lambda^{6}\displaystyle\iint_{(0,T)\times\omega_1} e^{-2s\alpha}\xi^{5}|\nabla^{N-1}\psi_1|^2 \right),
 \\ \ldots 
 \\  s^{2N+1}\lambda^{2N+2}\displaystyle\iint_{Q_T} e^{-2s\alpha}\xi^{2N+1}|\nabla\psi_1|^2 \\ \leqslant C\left(s^{2N-1}\lambda^{2N}\displaystyle\iint_{Q_T} e^{-2s\alpha}\xi^{2N-1}|\nabla^2\psi_1|^2\right . \\ \left . + s^{2N+1}\lambda^{2N+2}\displaystyle\iint_{(0,T)\times\omega_1} e^{-2s\alpha}\xi^{2N+1}|\nabla\psi_1|^2 \right).
  \end{array}\end{equation*}
We deduce by starting from the last inequality and using in cascade the other ones that 
\begin{equation}
\label{gauche1}
\begin{array}{l}
 s^{5}\lambda^{6}\displaystyle\iint_{Q_T} e^{-2s\alpha}\xi^{5}|\nabla^{N-1}\psi_1|^2\ldots +s^{2N+1}\lambda^{2N+2}\displaystyle\iint_{Q_T} e^{-2s\alpha}\xi^{2N+1}|\nabla\psi_1|^2 \\ 
 \leqslant C\left(s^{3}\lambda^{4}\displaystyle\iint_{Q_T} e^{-2s\alpha}\xi^{3}|\nabla^{N}\psi_1|^2  + s^{5}\lambda^{6}\displaystyle\iint_{(0,T)\times\omega_1} e^{-2s\alpha}\xi^{5}|\nabla^{N-1}\psi_1|^2\right.\\ \left .\ldots +s^{2N+1}\lambda^{2N+2}\displaystyle\iint_{(0,T)\times\omega_1}  e^{-2s\alpha}\xi^{2N+1}|\nabla\psi_1|^2 \right ).
 \end{array}\end{equation}

Combining \eqref{i1}, \eqref{gauche1} and using the definition of $I(s,\lambda,\phi)$ given in \eqref{defI}, we deduce that we can absorb the first term on the right-hand side of \eqref{gauche1} and obtain
\begin{equation*}
%\label{i2}
\begin{array}{l}
 s\lambda^{2}\displaystyle\iint_{Q_T} e^{-2s\alpha}\xi|\nabla^{N+1}\psi_1|^2\ldots +s^{2N+1}\lambda^{2N+2}\displaystyle\iint_{Q_T} e^{-2s\alpha}\xi^{2N+1}|\nabla\psi_1|^2  \\ 
\leqslant C\left( \lambda 
 \displaystyle\iint_{Q_T}e^{-2s\alpha^*}(s\xi_*)^{\frac{2 N (p+1)+5 p+3}{2 p}}
|\psi_1|^2   +\displaystyle\iint_{Q_T}e^{-2s\alpha}\sum_{i=1}^{N-1}|\nabla^i \psi_1|^2 \right .
\\ \left .  + s^{3}\lambda^{4}\displaystyle\iint_{(0,T)\times\omega_1} e^{-2s\alpha}\xi^{3}|\nabla^{N}\psi_1|^2+s^{5}\lambda^{6}\displaystyle\iint_{(0,T)\times\omega_1} e^{-2s\alpha}\xi^{5}|\nabla^{N}\psi_1|^2\right . \\ \left .\ldots +s^{2N+1}\lambda^{2N+2}\displaystyle\iint_{(0,T)\times\omega_1}  e^{-2s\alpha}\xi^{2N+1}|\nabla\psi_1|^2 \right ). \end{array}\end{equation*}
Absorbing the second term of the right-hand side, we deduce that for $s,\lambda$ large enough, we have
\begin{equation}
\label{i2}
\begin{array}{l}
 s\lambda^{2}\displaystyle\iint_{Q_T} e^{-2s\alpha}\xi|\nabla^{N+1}\psi_1|^2+\ldots +s^{2N+1}\lambda^{2N+2}\displaystyle\iint_{Q_T} e^{-2s\alpha}\xi^{2N+1}|\nabla\psi_1|^2  \\ 
\leqslant C\left( \lambda 
 \displaystyle\iint_{Q_T}e^{-2s\alpha^*}(s\xi_*)^{\frac{2 N (p+1)+5 p+3}{2 p}}
|\psi_1|^2  \right .
\\ \left  . + s^{3}\lambda^{4}\displaystyle\iint_{(0,T)\times\omega_1} e^{-2s\alpha}\xi^{3}|\nabla^{N}\psi_1|^2\right .\\ \left .\ldots +s^{2N+1}\lambda^{2N}\displaystyle\iint_{(0,T)\times\omega_1}  e^{-2s\alpha}\xi^{2N+1}|\nabla\psi_1|^2 \right ). \end{array}\end{equation}

{\bf Step 3: }
Now, we consider some open subset $\omega_{2}$ such that  $\omega_1\subset\subset \omega_{2}\subset\subset \omega_0$. We consider some function $\tilde\theta\in C^\infty(\Omega,\mathbb R)$ such that:
\begin{itemize}
\item $\textrm{Supp}(\tilde\theta)\subset \omega_2$,
\item $\tilde\theta=1$ on $\omega_1$,
\item $\tilde\theta\in [0,1]$.
\end{itemize}
Some integrations by parts give 
\begin{equation*}
\begin{array}{ll}
s^{3}\lambda^{4}\displaystyle\iint_{(0,T)\times\omega_1} e^{-2s\alpha}\xi^{3}|\nabla^{N}\psi_1|^2\leqslant s^{3}\lambda^{4}\displaystyle\iint_{(0,T)\times\omega_{2}}\theta e^{-2s\alpha}\xi^{3}|\nabla^{N}\psi_1|^2\\\displaystyle\leqslant Cs^3\lambda^4\iint_{(0,T)\times\omega_{2}}\left(|\nabla (\theta e^{-2s\alpha}\xi^{3})|.|\nabla^{N}\psi_1|.|\nabla^{N-1}\psi_1| \right . \\ \left .+|\theta e^{-2s\alpha}\xi^{3}|.|\nabla^{N+1}\psi_1|.|\nabla^{N-1}\psi_1|\right).
 \end{array}\end{equation*}
 From the definition of $\xi$ and $\alpha$ given in \eqref{defax} and \eqref{defax2}, we deduce that 
\begin{equation}\label{nth}|\nabla (\theta e^{-2s\alpha}\xi^{3})|\leqslant Cs\lambda e^{-2s\alpha}\xi^{4}.\end{equation}
 Combining this estimate with Young's inequality, we obtain that for any $\varepsilon>0$, there exists $C_\varepsilon>0$ such that for any $s$ and $\lambda$ large enough, we have
 \begin{equation}
 \label{azer}
\begin{array}{ll}
s^{3}\lambda^{4}\displaystyle\iint_{(0,T)\times\omega_1} e^{-2s\alpha}\xi^{3}|\nabla^{N}\psi_1|^2\leqslant C\left(\varepsilon s^{3}\lambda^{4}\displaystyle\iint_{(0,T)\times\omega_2} e^{-2s\alpha}\xi^{3}|\nabla^{N}\psi_1|^2 \right . \\\left .  + \varepsilon s\lambda^{2}\displaystyle\iint_{(0,T)\times\omega_2} e^{-2s\alpha}\xi|\nabla^{N+1}\psi_1|^2\right . \\\left .+C_\varepsilon s^{5}\lambda^{6}\displaystyle\iint_{(0,T)\times\omega_2} e^{-2s\alpha}\xi^{5}|\nabla^{N-1}\psi_1|^2 \right). \end{array}
\end{equation}
 Combining \eqref{i2} and \eqref{azer}, we can absorb the local terms in $|\nabla^{N+1}\psi_1|^2$ and $|\nabla^{N}\psi_1|^2$ to deduce  
 \begin{equation*}%\label{i21}
\begin{array}{c}
 s\lambda^{2}\displaystyle\iint_{Q_T} e^{-2s\alpha}\xi|\nabla^{N+1}\psi_1|^2\ldots +s^{2N+1}\lambda^{2N+2}\displaystyle\iint_{Q_T} e^{-2s\alpha}\xi^{2N+1}|\nabla\psi_1|^2  \\ 
\leqslant C\left( \lambda 
 \displaystyle\iint_{Q_T}e^{-2s\alpha^*}(s\xi_*)^{\frac{2 N (p+1)+5 p+3}{2 p}}
|\psi_1|^2  \right . \\\left . + s^{5}\lambda^{6}\displaystyle\iint_{(0,T)\times\omega_2} e^{-2s\alpha}\xi^{5}|\nabla^{N-1}\psi_1|^2\right .
\\ \left  .\ldots +s^{2N+1}\lambda^{2N+2}\displaystyle\iint_{(0,T)\times\omega_{2}}  e^{-2s\alpha}\xi^{2N+1}|\nabla\psi_1|^2 \right ). \end{array}\end{equation*}
We can perform exactly the same procedure on the terms 
 $$\begin{aligned} s^{5}\lambda^{6}\displaystyle\iint_{(0,T)\times\omega_2} e^{-2s\alpha}\xi^{5}|\nabla^{N-1}\psi_1|^2, \ldots, \\s^{2N-1}\lambda^{2N-2}\displaystyle\iint_{(0,T)\times\omega_2}  e^{-2s\alpha}\xi^{2N-1}|\nabla^2\psi_1|^2  \end{aligned}$$
 in order to obtain the following estimate:
  \begin{equation}
\label{i21}
\begin{array}{c}
 s\lambda^{2}\displaystyle\iint_{Q_T} e^{-2s\alpha}\xi|\nabla^{N+1}\psi_1|^2\ldots +s^{2N+1}\lambda^{2N+2}\displaystyle\iint_{Q_T} e^{-2s\alpha}\xi^{2N+1}|\nabla\psi_1|^2  \\ 
\leqslant C\left( \lambda 
 \displaystyle\iint_{Q_T}e^{-2s\alpha^*}(s\xi_*)^{\frac{2 N (r+1)+5 r+3}{2 r}}
|\psi_1|^2  \right.\hspace*{2cm}\\\hspace*{2cm}\left.+
s^{2N+1}\lambda^{2N+2}\displaystyle\iint_{(0,T)\times\omega_0}  e^{-2s\alpha}\xi^{2N+1}|\nabla\psi_1|^2 \right ). \end{array}
\end{equation}

{ \bf Step 4:}
 Since the weight  $(s\xi_*)^{2N-1}$ does not depend on the space variable, \textcolor{black}{$s\xi_*$ is bounded from below by a positive number}, and using the definition of $\alpha^*$ and $\xi_*$ given in \eqref{defaxs},
the following Poincar\'e's inequality holds: 

\begin{equation}\begin{array}{ll}\label{ine: poincarre psi2b}
 \lambda^{2N+2}\displaystyle\iint_{Q_T} e^{-2s\alpha^*}(s\xi_*)^{2N+1}|\psi_1|^2
 \\\leqslant  C\lambda^{2N+2}\displaystyle\iint_{Q_T} e^{-2s\alpha^*}(s\xi_*)^{2N+1}|\nabla\psi_1|^2\\
 \\\leqslant   C\lambda^{2N+2}\displaystyle\iint_{Q_T} e^{-2s\alpha}(s\xi)^{2N+1}|\nabla\psi_1|^2.
 \end{array}
\end{equation}
Combining \eqref{i21} and \eqref{ine: poincarre psi2b}, we deduce that 
 for $s$ large enough
\begin{equation}\label{yena}%\label{ine: neum carl proof2}
\begin{array}{l}
\lambda^{2}\displaystyle\iint_{Q_T} e^{-2s\alpha}s\xi|\nabla^{N+1}\psi_1|^2+\ldots +\lambda^{2N+2}\displaystyle\iint_{Q_T} e^{-2s\alpha}(s\xi)^{2N+1}|\nabla\psi_1|^2 \\+ \lambda^{2N+2}\displaystyle\iint_{Q_T}e^{-2s\alpha^*}(s\xi_*)^{2N+1}|\psi_1|^2
\\\leqslant C\left( \lambda 
 \displaystyle\iint_{Q_T}e^{-2s\alpha^*}(s\xi_*)^{\frac{2 N (p+1)+5 p+3}{2 p}}
|\psi_1|^2 \right . \\\left . +
\lambda^{2N+2}\displaystyle\iint_{(0,T)\times\omega_0}  e^{-2s\alpha}(s\xi)^{2N+1}|\nabla\psi_1|^2 \right ).

\end{array}
\end{equation}
We now fix $p\geqslant 2$ large enough such that 
$$\frac{2 N (p+1)+5 p+3}{2 p}<2N+1,$$
which is clearly possible since $\frac{2 N (p+1)+5 p+3}{2 p}\rightarrow N+\frac{5}{2}$ as $p\rightarrow \infty$ and $N\geqslant 3$ (so that $N+5/2<2N+1$).

Using  that $e^{-2s\alpha^*}(s\xi_*)^{\frac{2 N (p+1)+5 p+3}{2 p}}\leqslant Ce^{-2s\alpha}(s\xi)^{2N+1}$,  we deduce by absorbing the first  term of the right-hand side of \eqref{yena} that 

\begin{equation*}%\label{ine: neum carl proof2}
\begin{array}{l}
\lambda^{2}\displaystyle\iint_{Q_T} e^{-2s\alpha}(s\xi)|\nabla^{N+1}\psi_1|^2\ldots +\lambda^{2N+2}\displaystyle\iint_{Q_T} e^{-2s\alpha}(s\xi)^{2N+1}|\nabla\psi_1|^2\\ + \lambda^{2N+2}\displaystyle\iint_{Q_T}e^{-2s\alpha^*}(s\xi_*)^{2N+1}|\psi_1|^2\\
\leqslant C
\lambda^{2N+2}\displaystyle\iint_{(0,T)\times\omega_0}  e^{-2s\alpha}(s\xi)^{2N+1}|\nabla\psi_1|^2.
\end{array}
\end{equation*}
Going back to $\psi$ thanks to \eqref{sys:psi0}, we deduce \eqref{ine obs 1}.

\cqfd

\subsection{Algebraic resolubility}
\label{s:a}
In this section, we will derive a new Carleman inequality, adapted to the control problem with less controls we want to prove. 
%We assume here that $q\in\mathbb N^*$ (  in the following Lemma).

\begin{Lemme}\label{lemme:resol alg}
Let $m\in\mathbb N^*$ such that $m\leqslant d-1$.
Assume that the $\overline{u}$  is regular enough (for example of class $\mathcal C^\infty$).

Consider two partial differential operators $\mathcal L_1:C^\infty(\mathbb R^d)\rightarrow C^\infty(\mathbb R^d)^m$ and $\mathcal L_2:C^\infty(\mathbb R^d)\rightarrow C^\infty(\mathbb R^d)$ defined for every $\varphi\in \mathcal C^{\infty}(\mathbb R^d)$ by
\begin{equation*}%\label{def L*}
 \mathcal{L}_1\varphi:=B^*(\nabla\varphi)
 \mbox{ and }
 \mathcal{L}_2\varphi:=\partial_t\varphi+\Delta\varphi+(\overline{u}\cdot \nabla)\varphi.
\end{equation*}

Assume that \eqref{cdet} holds, \textcolor{black}{and let $q\in \mathbb N$ such that 
\begin{equation}\label{cdetq}\begin{aligned} &\textcolor{black}{\mbox{rank} }(\{B_1^* , ..., B_m^* \}\\& \cup\{((B^*\cdot\nabla)^\alpha \bar u_i(t, x))_{i\in \{1,\cdots,d\}}, \alpha\in\mathbb{N}^m,\,\alpha\not = 0,\, ||\alpha||_1\leqslant q  \})=d.\end{aligned}
\end{equation}}

There exists an open subset $(t_1,t_2)\times\widetilde\omega$ of $(0,T)\times\omega$ and there exist two partial differential operators $\mathcal M_1:C^\infty(\mathbb R^d)^m\rightarrow C^\infty(\mathbb R^d)^d$ (of order $1$ in time and $q+1$ in space) and $\mathcal M_2:C^\infty(\mathbb R)\rightarrow C^\infty(\mathbb R^d)^d$ (of order $0$ in time and $q$ in space) such that
\begin{equation}\label{M1M2}
\mathcal M_1\circ \mathcal L_1+\mathcal M_2\circ \mathcal L_2=\nabla\mbox{ in }\mathcal C^{\infty}({(t_1,t_2)\times\widetilde\omega}).
\end{equation}

\end{Lemme}

\textbf{Proof of Lemma \ref{lemme:resol alg}:}
\textcolor{black}{If $q=0$, necessarily, by condition \eqref{cdet}, we have $m=d$ and we can take $\mathcal M_1=(B^*)^{-1}$ and $\mathcal M_2=0$. 
We assume from now on that $q\in\mathbb N^*$.}
Let $j\in \{1,...,m\}$. We call $\mathcal L^j_1$ the $j-th$ line of $\mathcal L_1$.
We remark that
$$\begin{aligned}&(B^*_j\cdot\nabla)\mathcal L_2\varphi-(\partial_t+\Delta)\mathcal L^j_1\varphi-(\overline u\cdot\nabla)\mathcal L^j_1\varphi\\&=(B^*_j\cdot\nabla)(\overline{u}\cdot \nabla)\varphi-(\overline u\cdot\nabla)(B^*_j\cdot\nabla) \varphi
\\&=(\overline u\cdot \nabla)(B^*_j\cdot\nabla)\varphi+\sum_{k=1}^d ((B^*_j\cdot\nabla)\overline u_k) \partial_k \varphi\\&\hspace*{4cm}-(\overline u\cdot\nabla)(B^*_j\cdot\nabla) \varphi\\
&=\sum_{k=1}^d ((B^*_j\cdot\nabla)\overline u_k) \partial_k \varphi\\&=:\mathcal L_3^j.\end{aligned}$$

Now, for some $l\in \{1,...,m\}$, the same computations easily give

$$\begin{aligned}(B^*_l\cdot\nabla)\mathcal L_3^j\varphi-\sum_{k=1}^d ((B^*_j\cdot\nabla)\overline u_k) \partial_k \mathcal{L}^l_1\varphi
&=\sum_{k=1}^d ((B^*_l\cdot\nabla)(B^*_j\cdot\nabla)\overline u_k) \partial_k \varphi\\=:\mathcal L_4^{j,l}\varphi.\end{aligned}$$

Continuing this procedure, we can easily create two partial differential operators $\widetilde {\mathcal M}_1$ (of order $1$ in time and $q+1$ in space) and $\widetilde {\mathcal M}_2$ (of order $0$ in time and $q$ in space) such that 
\textcolor{black}{$$ \widetilde{\mathcal M_1} (\mathcal L_1(\varphi))(t_0,x_0)+\widetilde{\mathcal M_2}(\mathcal L_2(\varphi))(t_0,x_0)=\widetilde M(\overline{u})(\nabla \varphi)(t_0,x_0),$$
where $\widetilde M(\overline{u})(t_0,x_0)$ is a matrix composed by $d$ independent vectors of the family
$M(\overline{u})(t_0,x_0)$ with $\|\alpha\|_1\leqslant q$ (which is possible since \eqref{cdetq} is verified).  By continuity, there exists an open neighbourhood $(t_1,t_2)\times\widetilde\omega$ of $(t_0,x_0)$ in $(0,T)\times \omega$ and $C>0$ such that 
$|\det(\widetilde M(\overline{u}))|>C$ on $(t_1,t_2)\times\widetilde\omega$.
%, the space spanned by $M(\overline{u})$ has a maximal dimension on $(t_1,t_2)\times\widetilde\omega$, so that a matrix $\widetilde{M}(\overline{u})$ composed by some of free subfamily  admits a left inverse at any point of on $(t_1,t_2)\times\widetilde\omega$. 
We call $\widetilde{M}(\overline{u})^{-1}(t,x)$ the inverse of  $\widetilde M(\overline{u})(t,x)$ for $(t,x)\in (t_1,t_2)\times\widetilde\omega$. 
Then,  is is clear that  $\mathcal M_1:=\widetilde{M}(\overline{u})^{-1}\widetilde {\mathcal M}_1$ 
and $\mathcal M_2:= \widetilde{M}(\overline{u})^{-1}\widetilde {\mathcal M}_2$ verify
 \eqref{M1M2} and have $C^\infty$ coefficients on $(t_1,t_2)\times\widetilde\omega$.
}
%where $M_{q}$ is defined in \eqref{mtm}.
%Under condition \eqref{cdet}, $M_{q}$ is of maximal rank on $(t_1,t_2)\times\widetilde\omega$, so that it admits a left inverse at any point of on $(t_1,t_2)\times\widetilde\omega$. 
%We call $ M_q(\overline{u})^{-1}$ any of its left inverses. Then,  is is clear that  $\mathcal M_1:=M_{q}^{-1}\widetilde {\mathcal M}_1$ and $\mathcal M_2:=M_{q}^{-1}\widetilde {\mathcal M}_2$ verify \eqref{M1M2} and have $C^\infty$ coefficients on $(t_1,t_2)\times\widetilde\omega$.
\cqfd

We now have all the tools to deduce our final Carleman inequality:
\begin{prop}\label{prop ine obs 2}
Assume that Condition \eqref{cdet} 
and the hypotheses of Proposition \ref{prop ine obs 1} hold.
Then, for all $\eta\in(0,1)$, 
there exists $p\geqslant 2$, $C>0$ and $K>0$ such that for every 
  $\psi^0\in L^2(\Omega)$, 
 the corresponding solution $\psi$ to System \eqref{sys:dual lin} satisfies
 \begin{equation}\label{ine obs 2}
 \begin{array}{c}
 \displaystyle\int_{\Omega}\psi(0)^2dx
 +\displaystyle\iint_{Q_T} e^{\frac{-2K}{\eta(T-t)^p}}
\{\psi^2+|\partial_t \psi|^2+\ldots +|\partial^{\lfloor \frac{N+1}{2}\rfloor}_{t\ldots t} \psi|^2\\+|\nabla\psi|^2+\ldots +|\nabla^{N+1}\psi|^2\}
\\\leqslant Ce^{K/T^p} \displaystyle\iint_{(0,T)\times \omega_0}e^{\frac{-2K}{(T-t)^p}}|B^*(\nabla \psi)|^2.
\end{array} 
\end{equation} 
\end{prop}
\textbf{Proof of Proposition \ref{prop ine obs 2}.}
%We assume that $q\in\mathbb N^*$ \textcolor{black}{ is such that \eqref{cdetq} holds} (\textcolor{black}{the case $q=0$ is trivial as already explained}).
Let $\omega_1$ some open subset strongly included in $\omega_0$.
Combining Proposition \ref{prop ine obs 1}, Lemma \ref{lemme:resol alg} (that is still true by replacing $\omega_0$ by $\omega_1$), and the fact that  any solution $\psi$ of \eqref{sys:dual lin} verifies by definition $\mathcal L_2 \psi=0$, we deduce that,
for any  $\psi^0\in L^2(\Omega)$, 
 the corresponding solution $\psi$ to System \eqref{sys:dual lin} satisfies
 \begin{equation*}%\label{ine: neum carl proof2}
\begin{array}{c}
\lambda^{2}\displaystyle\iint_{Q_T} e^{-2s\alpha-2\mu s\alpha^*}(s\xi)|\nabla^{N+1}\psi|^2+\ldots \\+\lambda^{2N+2}\displaystyle\iint_{Q_T} e^{-2s\alpha-2\mu s\alpha^*}(s\xi)^{2N+1}|\nabla\psi|^2\\ + \lambda^{2N+2}\displaystyle\iint_{Q_T}e^{-2s\alpha^*-2\mu s\alpha^*}(s\xi_*)^{2N+1}|\psi|^2\\
\leqslant C
\lambda^{2N+2}\displaystyle\iint_{Q_T}\tilde\theta  e^{-2s\alpha-2\mu s\alpha^*}\left(s\xi\right)^{2N+1}|\mathcal{M}_1B^*(\nabla \psi)|^2,
\end{array}
\end{equation*}
where $\mathcal M_1$ is a linear partial differential operator of order $1$ in time and $q+1$ in space, and $\tilde\theta\in C^\infty(\Omega,\mathbb R)$ such that:
\begin{itemize}
\item $\tilde\theta=1$ on $\omega_1$,
\item $\textrm{Supp}(\tilde\theta)\subset \omega_0$,
\item $\tilde\theta\in [0,1]$.
\end{itemize}
We first remark that 
 \begin{equation*}%\label{ine: neum carl proof2}
\begin{array}{c}
\displaystyle\iint_{Q_T}\theta  e^{-2s\alpha-2\mu s\alpha^*}\left(s\xi\right)^{2N+1}|\mathcal{M}_1B^*(\nabla \psi)|^2\\
\leqslant C\displaystyle\iint_{Q_T}\theta  e^{-2s\alpha-2\mu s\alpha^*}\left(s\xi\right)^{2N+1}\left                                                                                                                                                                                                                     (\sum_{i=0}^{q+1}\left (|\nabla^iB^*\nabla \psi|^2+|\partial_t\nabla^iB^*\nabla \psi|^2 \right) \right).
\end{array}
\end{equation*}
Using that  $\psi$ verifies \eqref{sys:dual lin}, we can deduce that 
 \begin{equation*}%\label{ine: neum carl proof2}
\begin{array}{c}
\lambda^{2N+2}\displaystyle\iint_{Q_T}\theta  e^{-2s\alpha-2\mu s\alpha^*}\left(s\xi\right)^{2N+1}|\mathcal{M}_1B^*(\nabla \psi)|^2\\
\leqslant C\lambda^{2N+2}\displaystyle\iint_{Q_T}\theta  e^{-2s\alpha-2\mu s\alpha^*}\left(s\xi\right)^{2N+1}\left                                                                                                                                                                                                                     (\sum_{i=0}^{q+3}|\nabla^iB^*\nabla \psi|^2 \right).
\end{array}
\end{equation*}
Some integrations by parts give 
 \begin{equation*}%\label{ine: neum carl proof2}
\begin{array}{c}
\lambda^{2N+2}\displaystyle\iint_{Q_T}\tilde\theta  e^{-2s\alpha-2\mu s\alpha^*}\left(s\xi\right)^{2N+1}|\nabla B^*(\nabla \psi)|^2\\
\leqslant C
\lambda^{2N+2}\displaystyle\iint_{Q_T}\tilde\theta  e^{-2s\alpha-2\mu s\alpha^*}\left(s\xi\right)^{2N+1}|B^*(\nabla \psi)||\nabla^3 \psi|\\
+C\lambda^{2N+2}\displaystyle\iint_{Q_T}|\nabla(\tilde\theta  e^{-2s\alpha-2\mu s\alpha^*}\left(s\xi\right)^{2N+1})||B^*(\nabla \psi)||\nabla^2 \psi|.
\end{array}
\end{equation*}
Let $\varepsilon>0$. Young's inequality gives
 \begin{equation*}%\label{ine: neum carl proof2}
\begin{array}{c}
\lambda^{2N+2}\displaystyle\iint_{Q_T}\tilde\theta  e^{-2s\alpha-2\mu s\alpha^*}\left(s\xi\right)^{2N+1}|B^*(\nabla \psi)||\nabla^3 \psi|\\
\leqslant C_{\varepsilon}\lambda^{2N+6}\displaystyle\iint_{(0,T)\times\omega_0}  e^{-2s\alpha-2\mu s\alpha^*}\left(s\xi\right)^{2N+5}|B^*(\nabla \psi)|^2
\\+\varepsilon\lambda^{2N-2}\displaystyle\iint_{(0,T)\times\omega_0}  e^{-2s\alpha-2\mu s\alpha^*}\left(s\xi\right)^{2N-3}|\nabla^3 \psi|^2
\end{array}
\end{equation*}
and also, by \eqref{nth},
 \begin{equation*}%\label{ine: neum carl proof2}
\begin{array}{c}
\lambda^{2N+2}\displaystyle\iint_{Q_T}|\nabla(\tilde\theta  e^{-2s\alpha-2\mu s\alpha^*}\left(s\xi\right)^{2N+1})||B^*(\nabla \psi)||\nabla^2 \psi|\\
\leqslant C\lambda^{2N+3}\displaystyle\iint_{Q_T}\tilde\theta  e^{-2s\alpha-2\mu s\alpha^*}\left(s\xi\right)^{2N+2}|B^*(\nabla \psi)||\nabla^2 \psi|\\
\leqslant C_{\varepsilon}\lambda^{2N+6}\displaystyle\iint_{(0,T)\times\omega_0} e^{-2s\alpha-2\mu s\alpha^*}\left(s\xi\right)^{2N+5}\textcolor{black}{|B^*(\nabla \psi)|^2}
\\+\varepsilon\lambda^{2N}\displaystyle\iint_{(0,T)\times\omega_0}  e^{-2s\alpha-2\mu s\alpha^*}\left(s\xi\right)^{2N-1}\textcolor{black}{|\nabla^2 \psi|^2}.
\end{array}
\end{equation*} 
Thus, by taking $\varepsilon$ small enough, we deduce that
 \begin{equation*}%\label{ine: neum carl proof2}
\begin{array}{c}
\lambda^{2}\displaystyle\iint_{Q_T} e^{-2s\alpha-2\mu s\alpha^*}(s\xi)|\nabla^{N+1}\psi|^2+\ldots \\+\lambda^{2N+2}\displaystyle\iint_{Q_T} e^{-2s\alpha-2\mu s\alpha^*}(s\xi)^{2N+1}|\nabla\psi|^2\\ + \lambda^{2N+2}\displaystyle\iint_{Q_T}e^{-2s\alpha^*-2\mu s\alpha^*}(s\xi_*)^{2N+1}|\psi|^2\\
\leqslant C
\lambda^{2N+6}\displaystyle\iint_{(0,T)\times\omega_0}  e^{-2s\alpha-2\mu s\alpha^*}\left(s\xi\right)^{2N+5}|B^*(\nabla \psi)|^2\\
+C\lambda^{2N+2}\displaystyle\iint_{Q_T}\tilde\theta  e^{-2s\alpha-2\mu s\alpha^*}\left(s\xi\right)^{2N+1}\left(\sum_{i=2}^{q+3}|\nabla^iB^*\nabla \psi|^2\right).
\end{array}
\end{equation*}
By iterating this process for $i=2,\ldots, q+3$, we can get rid of the sum in the right-hand side and obtain
 \begin{equation*}%\label{ine: neum carl proof2}
\begin{array}{c}
\lambda^{2}\displaystyle\iint_{Q_T} e^{-2s\alpha-2\mu s\alpha^*}(s\xi)|\nabla^{N+1}\psi|^2+\ldots \\+\lambda^{2N+2}\displaystyle\iint_{Q_T} e^{-2s\alpha-2\mu s\alpha^*}(s\xi)^{2N+1}|\nabla\psi|^2\\ + \lambda^{2N+2}\displaystyle\iint_{Q_T}e^{-2s\alpha^*-2\mu s\alpha^*}(s\xi_*)^{2N+1}|\psi|^2\\
\leqslant C 
\lambda^{2N+2+4(q+2)}\displaystyle\iint_{(0,T)\times\omega_0}  e^{-2s\alpha-2\mu s\alpha^*}\left(s\xi\right)^{2N+1+4(q+2)}|B^*(\nabla \psi)|^2 .

\end{array}
\end{equation*}
\textcolor{black}{We deduce that 
 \begin{equation*}%\label{ine: neum carl proof2}
\begin{array}{c}
\lambda^{2}\displaystyle\iint_{Q_T} e^{-2(1+\mu) s\alpha^*}(s\xi_*)|\nabla^{N+1}\psi|^2+\ldots\\ +\lambda^{2N+2}\displaystyle\iint_{Q_T} e^{-2(1+\mu) s\alpha^*}(s\xi_*)^{2N+1}|\nabla\psi|^2\\ + \lambda^{2N+2}\displaystyle\iint_{Q_T}e^{-2(1+\mu)\mu s\alpha^*}(s\xi_*)^{2N+1}|\psi|^2\\
\leqslant C 
\lambda^{2N+2+4(q+2)}\displaystyle\iint_{(0,T)\times\omega_0}  e^{-2\mu s\alpha^*}\left(s\xi^*\right)^{2N+1+4(q+2)}|B^*(\nabla \psi)|^2,
\end{array}
\end{equation*}
where $\xi^*=\max\limits_{\Omega}\xi$. 
Defining 
$$\begin{aligned}\tilde\alpha^*=\left\{\begin{array}{ll}
\alpha^*(T/2)&\mbox{ on }(0,T/2),\\
\alpha^*&\mbox{ on }(T/2,T), 
 \end{array}\right. \\
 \tilde\xi_*=\left\{\begin{array}{ll}
\xi_*(T/2)&\mbox{ on }(0,T/2),\\
\xi_*&\mbox{ on }(T/2,T), 
 \end{array}\right. \\
 \tilde\xi^*=\left\{\begin{array}{ll}
\xi^*(T/2)&\mbox{ on }(0,T/2),\\
\xi^*&\mbox{ on }(T/2,T), 
 \end{array}\right.\end{aligned}$$
 then, for $s$ and $\lambda$ large enough, using usual energy estimates, 
  \begin{equation*}%\label{ine: neum carl proof2}
\begin{array}{c}
\lambda^{2}\displaystyle\iint_{Q_T} e^{-2(1+\mu) s\tilde\alpha^*}(s\tilde\xi_*)|\nabla^{N+1}\psi|^2+\ldots \\+\lambda^{2N+2}\displaystyle\iint_{Q_T} e^{-2(1+\mu) s\tilde\alpha^*}(s\tilde\xi_*)^{2N+1}|\nabla\psi|^2\\ + \tilde\lambda^{2N+2}\displaystyle\iint_{Q_T}e^{-2(1+\mu)\mu s\tilde\alpha^*}(s\tilde\xi_*)^{2N+1}|\psi|^2\\
\leqslant C 
\lambda^{2N+2+4(q+2)}\displaystyle\iint_{(0,T)\times\omega_0}  e^{-2\mu s\tilde\alpha^*}\left(s\tilde\xi^*\right)^{2N+1+4(q+2)}|B^*(\nabla \psi)|^2 .
\end{array}
\end{equation*}
Fixing $s$ and $\lambda$, using \eqref{defax} and \eqref{defax2}, and remarking that $\tilde{\xi}^*$ does not depend on $\mu$, we deduce that there exists $R>0$ such that for any $\mu>0$ large enough,
 \begin{equation*}%\label{ine: neum carl proof2}
\displaystyle\iint_{Q_T} e^{-\frac{2(2+\mu) R}{(T-t)^p}}\{|\nabla^{N+1}\psi|^2+\ldots +|\psi|^2\}
\leqslant C
\displaystyle\iint_{(0,T)\times\omega_0}  e^{-\frac{2(\mu-1) R}{(T-t)^p}}|B^*(\nabla \psi)|^2 .
\end{equation*}
We remark that the fact that $\psi$ verifies \eqref{sys:dual lin} enables us to add all the derivatives in time on the left-hand side. Hence, we can conclude by fixing $\eta\in(0,1)$, introducing $K=(\mu-1) R$ (for $\mu>1$), and taking $\mu>0$ large enough so that 
$$\frac{(\mu-1) R}{\eta}>(2+\mu) R,
$$
which is always possible since the ratio  $(\mu-1)/(2+\mu)$ tends to $1$ as $\mu\rightarrow \infty$.
 }
%Inequality \eqref{ine obs 2} is easily deduced by replacing the space-dependent weights by their infimum in space in the left-hand-side and their supremum in the right-hand side, fixing $s$ and $\lambda$ large enough, then choosing $\mu$ large enough  (depending on $||\eta_0||_{\infty(\bar\Omega)}$) with respect to the parameter $\eta\in(0,1)$, applying usual energy estimates and remarking that the fact that $\psi$ verifies \eqref{sys:dual lin} enables us to add all the derivatives in time on the left-hand side.

\cqfd

\subsection{Regular control}

\label{s:rc}

Our goal in this section is to construct regular enough controls. Remind that $\theta$ is defined in
\eqref{deftheta}.

\begin{prop}\label{prop cont reg} 
Let $r\in\mathbb N$.
Assume that Condition \eqref{cdet} holds.

Under the hypotheses of Proposition \ref{prop ine obs 1}, 
System 
\begin{equation}\label{syst lin primal 2}
\left\{\begin{array}{ll}
\partial_ty=\Delta y+\Div(\overline u y)+\Div(\theta Bv)&\mbox{in } Q_T,\\
y=0&\mbox{on } \Sigma_T,\\
y(0,\cdot)=y^0&\mbox{in }\Omega,
\end{array}\right.
\end{equation} 
 is null controllable at time $T$, 
\textit{i.e.} for every $y^0\in L^2(\Omega)$, 
there exists a control $v\in L^2(Q_T)^m$ such that the solution $z$ 
to System \eqref{syst lin primal 2} satisfies $z(T)\equiv 0$ in $\Omega$. 
Moreover, we can choose $\textcolor{black}{v\in} L^2((0,T)$, $H^{2r+2}(\Omega))^m\cap H^{r+1}((0,T),L^{2}(\Omega))^m$ with
\begin{equation*}%\label{exp control}
\|v\|_{L^2((0,T),H^{2r+2}(\Omega))^m\cap H^{r+2}((0,T),L^{2}(\Omega))^m}
\leqslant Ce^{K/T^p}\|y^0\|_{L^2(\Omega)},
\end{equation*}
where $K$ is the constant in \eqref{ine obs 2}.
\end{prop}

\textbf{Proof of Proposition \ref{prop cont reg}.}
%Since $\overline u_i$ is regular and with spatial support strongly included in $\Omega$, without loss of generality, we may assume that $y^0\in H^{2r+1}(\Omega)$  with  
%$y^0=\Delta y^0\ldots = \Delta^ry^0=0$  on $\partial\Omega$ (by letting the system evolve freely during a short time).
Let $k\in\mathbb{N}^*$ and let us consider the following optimal control problem\begin{equation}\label{probleme minimisation}
 \left\{\begin{array}{l}
         \mathrm{minimize}~ J_k(v):=\dfrac12\|\widetilde{\rho}^{-1/2} v\|^2_{L^2(Q_T)^m}
 +\dfrac{k}{2}\displaystyle\int_{\Omega}| z(T)|^2dx,\\
 v\in \mathcal{U}:=%L^2(Q_T,\widetilde{\rho}^{-1})^m =
 \{w\in L^2(Q_T)^m:\widetilde{\rho}^{-1/2} w\in L^2(Q_T)^m\},
        \end{array}
\right.
\end{equation}
where  $\widetilde{\rho}:=e^{\frac{-2K}{(T-t)^p}}$ (for the $K>0$ given by Proposition \ref{prop ine obs 2} with $N$ an even number to be chosen later and some fixed $\eta \in(1/2,1)$)
 and $z$ is the solution in $W(0,T)$ to 
\begin{equation*}
\left\{\begin{array}{ll}
\partial_tz=\mathcal{A} z+\mathcal{B}v&\mbox{in } Q_T,\\
y=0&\mbox{on } \Sigma_T,\\
y(0,\cdot)=y^0&\mbox{in }\Omega,
\end{array}\right.
\end{equation*}
where
\begin{equation}\label{syst lin primal reg}
\left\{\begin{array}{l}
\mathcal{A}:=\Delta +\Div(\overline u ~\cdot~),\\
\mathcal{B}:=\Div(B\theta ~\cdot~).
\end{array}\right.
\end{equation}
Here, $\mathcal U$ is endowed with its natural weighted $L^2$-norm.

%The multiplication by $\widetilde{\rho}_k$ is an bounded operator.
%Moreover, t
The functional $J_k:\mathcal{U}\rightarrow\mathbb{R}^+$ is differentiable, coercive and strictly 
convex on the space $\mathcal{U}$. 
Therefore, following \cite[[p. 116]{lionscontrole},
there exists a unique solution 
to the optimal control problem \eqref{probleme minimisation} 
and the optimal control $v_k$ 
is characterized thanks to the solution $z_k$ of the primal system by 
\begin{equation}\label{srt}
 \left\{\begin{array}{ll}
       \partial_tz_k=\mathcal{A}z_k+\mathcal{B}v_k&\mathrm{in}~ Q_T,\\
       z_k=0&\mathrm{on}~\Sigma_T,\\
       z_k(0,\cdot)=y^0&\mathrm{in}~\Omega,
        \end{array}
\right.
\end{equation}
the solution $\varphi_k$ to the dual system
\begin{equation}\label{system dual fursikov}
\left\{\begin{array}{ll}
         -\partial_t\varphi_k=\mathcal{A}^*\varphi_k&\mathrm{in}~ Q_T,\\
                \varphi_k=0&\mathrm{on}~\Sigma_T,\\
       \varphi_k(T,\cdot)=k  z_k(T,\cdot)&\mathrm{in}~\Omega
        \end{array}
\right.
\end{equation}
and the relation
\begin{equation}\label{carac controle}
\left\{\begin{array}{l}
        v_k=-\widetilde{\rho}
        \mathcal{B}^*\varphi_k\mathrm{~in~}Q_T,\\
       v_k\in \mathcal{U}.
        \end{array}
\right.
\end{equation}

The characterization  \eqref{srt}, \eqref{system dual fursikov} and \eqref{carac controle} of the minimizer $v_k$ of $J_k$  in $\mathcal{U}$ 
leads to the following computations:
\begin{equation}\label{estim J_k1}
\begin{array}{rcl}
 J_k(v_k)&=&-\dfrac12\langle 
 \mathcal{B}^*\varphi_k,v_k\rangle_{L^2(Q_T)^m}
 +\dfrac12\langle z_k(T),\varphi_{k}(T)\rangle_{L^2(\Omega)}\vspace*{3mm}\\
 &=&-\dfrac12\displaystyle\int_{0}^T\langle \varphi_k,\mathcal{B}v_k\rangle_{L^2(\Omega)}
  +\dfrac12\displaystyle\int_{0}^T\{\langle z_k,\partial_t\varphi_k\rangle_{L^2(\Omega)}\\
 &&~~~~~~~~~~~~~~~~~~~~~~~~~~~~~~~~~~~~~~~~~~~~~~~
  +\langle\partial_tz_k,\varphi_k\rangle_{L^2(\Omega)}\}
   +\dfrac12\langle y^0,\varphi_{k}(0,\cdot)\rangle_{L^2(\Omega)}\\
 &=& \dfrac12\langle y^0,\varphi_{k}(0,\cdot)\rangle_{L^2(\Omega)}.
 \end{array}
 \end{equation}

Moreover, using \eqref{ine obs 2} %with $N=2s$ 
and the expression of $\widetilde{\rho}$, 
%inequality $$\widetilde{\rho}_k(t)\leqslant e^{\frac{-2K}{(T-t)^p}},$$ 
we infer 
  \begin{equation}\label{pasutilise}
  \begin{array}{rcl}
 \|\varphi_k(0,\cdot)\|_{L^2(\Omega)}
&  \leqslant &
Ce^{K/T^p} \|\widetilde{\rho}^{-1/2}v_k\|_{L^2(Q_T)^m}.
\end{array}
 \end{equation}

Now, using the definition of $J_k$, 
the expression \eqref{estim J_k1}, the inequality \eqref{pasutilise} and the Cauchy-Schwartz inequality, we infer 
  \begin{equation*}
  \begin{array}{cl}
 \|\varphi_k(0,\cdot)\|_{L^2(\Omega)}^2
 &\leqslant Ce^{2K/T^p}J_k(v_k)\leqslant Ce^{2K/T^p} \|\varphi_k(0,\cdot)\|_{L^2(\Omega)}\|y^0\|_{L^2(\Omega)},
\end{array}
 \end{equation*}
from which we deduce 
  \begin{equation}\label{estim J_k2}
 \|\varphi_k(0,\cdot)\|_{L^2(\Omega)}\leqslant Ce^{2K/T^p}\|y^0\|_{L^2(\Omega)}.
 \end{equation}
 Then, using  \eqref{estim J_k1} and \eqref{estim J_k2}, we deduce 

\begin{equation}\label{J_k(V) borne}
 J_k(v_k)\leqslant{Ce^{2K/T^p}}\|y^0\|_{L^2(\Omega)}^2.
\end{equation}
Furthermore, we have (see \cite[p. 116]{lionscontrole})
\begin{equation}\label{norme z bornee}
\begin{array}{rl}
 \|z_k\|_{W(0,T)}
 &\leqslant C\left(\| \mathcal{B}v_k\|_{L^2((0,T),H^{-1}(\Omega))}+\|y^0\|_{L^2(\Omega)}\right),\\
  &\leqslant C\left(\|\widetilde{\rho}^{-1/2}v_k\|_{L^2(Q_T)^m}+\|y^0\|_{L^2(\Omega)}\right),\\
  &\leqslant C(1+Ce^{K/T^p})  \|y^0\|_{L^2(\Omega)},
\end{array}
\end{equation}
where $C$ does not depend on $y^0$ and $k$. 
Then, using  inequalities \eqref{J_k(V) borne} and \eqref{norme z bornee}, we deduce 
that there exist subsequences, which are still denoted $v_k$, $z_k$, 
such that the following weak convergences hold:
\begin{equation*}\left\{
 \begin{array}{ll}
  v_k\rightharpoonup v~&\mathrm{in}~\mathcal{U},\\
  z_k\rightharpoonup z~&\mathrm{in}~W(0,T),\\
    z_k(T)\rightharpoonup   0~&\mathrm{in}~L^2(\Omega).\\
 \end{array}
\right.
\end{equation*}
Passing to the limit in $k$, $z$ is solution to System \eqref{syst lin primal reg}. 
Moreover, using the expression of $J_k$ given in \eqref{probleme minimisation} and inequality \eqref{J_k(V) borne}, we deduce by letting $k$ going to $\infty$ that $ z(T)\equiv0$ in $\Omega$. 
Thus the solution $z$ to System \eqref{syst lin primal reg} with control $v\in \mathcal{U}$ 
satisfies $z(T)\equiv0$ in $\Omega$ 
and using (\ref{J_k(V) borne}), we obtain the inequality 
\begin{equation*}
%\label{normu}
 \|v\|_{\mathcal{U}}^2\leqslant Ce^{2K/T^p} \|y^0\|_{L^2(\Omega)}^2.
\end{equation*}
Since $\tilde\rho^{-1}\geqslant 1$ , using the definition of the norm on $\mathcal U$, we also deduce that 
\begin{equation*}
 \|v\|_{L^2(Q_T)^m}^2\leqslant Ce^{2K/T^p} \|y^0\|_{L^2(\Omega)}^2.
\end{equation*}

%\textcolor{black}{Je ne comprends rien a partir d'ici: on a pas le Carleman avec les poids $\widetilde{\rho}_k$. En fait, a mon avis, ce qu il faut faire c est utiliser la version sans k de \eqref{pasutilise}
%et en deduire les regularites voulus en disant que dans \eqref{ine obs 2} si le lrh est fini alors le lhs est fini.
Now, let us explain why the controls are more regular. First of all, using the fact that $\varphi_k$ verifies \eqref{system dual fursikov}, we deduce that 
\textcolor{black}{for any $j\in\mathbb{N}$,
%$$||\mathcal B^*\varphi_k||^2_{L^2(Q_T)}\leqslant C\|\partial_t \varphi_k\|^2_{L^2(Q_T)}.$$
%$$||\mathcal B^*\varphi_k||^2_{L^2(Q_T)}\leqslant C\|\partial_t \varphi_k\|^2_{L^2(Q_T)}.$$
$$\|\mathcal B^* \partial_t^j \varphi_k (t, .)\|^2_{ L^2 (\Omega)}\leqslant  C\|\partial_t^{j+1}\varphi_k (t, .)\|^2_{ L^2 (\Omega)},\forall t \in (0, T ).$$}
%
%Hence,  for each $i\in\{1,...,\frac{N}{2}-1\}$ and $k\in\mathbb{N}$, using inequalities similar to \eqref{rhot} and \eqref{rhot2}, we deduce that for any $\varepsilon>0$, there exists $C>0$ such that 
%%and the fact that for any $j\in\mathbb N$, 
%%$$\|\partial_t^j(\widetilde{\rho}_k)\|\leqslant C\|\partial_t^j(\widetilde{\rho}_k)\|
% %(with $C$ that might depend on $j$, $K$ and $T$ but not on $k$)
%\begin{equation}\label{cre}\begin{aligned}
% \|\partial_t^iv_k\|_{L^2(Q_T)^m}^2&=\iint_{Q_T}\partial_t^i\left(|-\tilde\rho \mathcal B^*\varphi_k \right|)^2\\&\leqslant C\iint_{Q_T}\tilde \rho^{2-2\varepsilon}\textcolor{black}{\{|\partial_t^{i-1} \varphi_k |^2+|\partial_t^{i} \varphi_k |^2\}}
%\\&\leqslant C\iint_{Q_T}\tilde \rho^{2-2\varepsilon-\frac{1}{\eta}}\tilde \rho^{\frac{1}{\eta}}\textcolor{black}{\{|\partial_t^{i-1} \varphi_k |^2+|\partial_t^{i} \varphi_k |^2\}}.
%\end{aligned}
%\end{equation}
Hence,  for each $i\in\{1,...,\frac{N}{2}-1\}$ and $k\in\mathbb{N}$, using inequalities similar to \eqref{rhot} and \eqref{rhot2}, we deduce that for any $\varepsilon>0$, there exists $C>0$ such that 
%and the fact that for any $j\in\mathbb N$, 
%$$\|\partial_t^j(\widetilde{\rho}_k)\|\leqslant C\|\partial_t^j(\widetilde{\rho}_k)\|
 %(with $C$ that might depend on $j$, $K$ and $T$ but not on $k$)
\begin{equation}\label{cre}\begin{aligned}
 \|\partial_t^iv_k\|_{L^2(Q_T)^m}^2&=\iint_{Q_T}\partial_t^i\left(|-\tilde\rho \mathcal B^*\varphi_k \right|)^2\\&\leqslant C\iint_{Q_T}\tilde \rho^{2-2\varepsilon}\textcolor{black}{\{|\varphi_k|^2+\cdots+|\partial_t^{i+1} \varphi_k|^2\}}
\\&\leqslant C\iint_{Q_T}\tilde \rho^{2-2\varepsilon-\frac{1}{\eta}}\tilde \rho^{\frac{1}{\eta}}\textcolor{black}{\{|\varphi_k|^2+\cdots+|\partial_t^{i+1} \varphi_k|^2\}}.
\end{aligned}
\end{equation}
Now, we fix $\varepsilon>0$ small enough (with respect to $\eta$) such that 
$2-2\varepsilon-\frac{1}{\eta}\geqslant 0$.
With this choice of $\varepsilon$, we infer that 
$\tilde \rho^{2-2\varepsilon-\frac{1}{\eta}}\leqslant 1$.
Hence, using \eqref{cre} together with \eqref{ine obs 2} and \eqref{J_k(V) borne}, we deduce that, for each $i\in\{0,...,\frac{N}{2}-1\}$, $\|\partial_t^{i}v_k\| \in L^2(Q_T)$ and
\begin{equation*}%\label{cre}
\begin{aligned}
 \|\partial_t^iv_k\|_{L^2(Q_T)^m}^2&\leqslant  C\iint_{Q_T} e^{\frac{-2K}{\eta(T-t)^p}}\textcolor{black}{\{|\varphi_k|^2+\cdots+|\partial_t^{i+1} \varphi_k|^2\}}\\ &\leqslant C\iint_{Q_T}e^{\frac{-2K}{(T-t)^p}}|\theta \mathcal B^*(\varphi_k)|^2\\&\leqslant C \|v_k\|_{\mathcal{U}}^2\\&\leqslant Ce^{2K/T^p} \|y^0\|_{L^2(\Omega)}^2.
\end{aligned}
\end{equation*}

$$\|\mathcal B^* \partial_t^j \phi_k (t, .)\|^2_{ L^2 (\Omega)}\leqslant  C\|\partial_t^{j+1}\phi_k (t, .)\|^2_{ L^2 (\Omega)},\forall t \in (0, T ),$$

Thus, extracting one more time a subsequence if necessary and letting $k$ go to $+\infty$, we deduce that for each $i\in\{1,...,\frac{N}{2}-1\}$,
\[
 \|\partial_t^iv\|_{L^2(Q_T)^m}
\leqslant Ce^{2K/T^p} \|y^0\|_{L^2(\Omega)}^2.
\] 
We similarly deduce that, for each $i\in\{1,...,N-2\}$,
\[
 \|\nabla^iv\|_{L^2(Q_T)^{m\times i\times d}}
\leqslant Ce^{2K/T^p} \|y^0\|_{L^2(\Omega)}^2.
\] 

The proof is completed by setting $r=\frac{N}{2}+1$.
\cqfd

\section{Controllability to the trajectories}
\label{s:nl}
Let $r\in\mathbb N$. We use the strategy developed in \cite{MR3023058}, modifying it slightly to fit our case. Usual interpolation estimates (see \cite[Section 13.2,
p. 96]{lions1968problemes}) show that
$$\begin{aligned}L^2((0,T),H^{2r+2}(\Omega))\cap H^{r+1}((0,T),L^{2}(\Omega))\\ \hookrightarrow L^2((0,T),H^{2r+2}(\Omega))\cap H^{1}((0,T),H^{2r}(\Omega)),\end{aligned}$$
from which we deduce 
$$L^2((0,T),H^{2r+2}(\Omega))\cap H^{r+1}((0,T),L^{2}(\Omega)) \hookrightarrow L^\infty((0,T),H^{2r}(\Omega)).$$
Now, there exists $R>0$ large enough such that by Sobolev embeddings, we have 
$$L^2((0,T),H^{2R+2}(\Omega))\cap H^{R+1}((0,T),L^{2}(\Omega)) \hookrightarrow L^\infty((0,T),W^{1,\infty}(\Omega)).$$

Hence, from Proposition \ref{prop cont reg} and Remark \ref{eqt}, for any $y^0\in L^2(\Omega)$, there exists a control $v\in L^\infty((0,T),W^{1,\infty}(\Omega))^m$ such that the solution $y$ 
to System \eqref{reale} satisfies $y(T)\equiv 0$ in $\Omega$ and 
\begin{equation*}%\label{exp control2}
\|v\|_{ L^\infty((0,T),W^{1,\infty}(\Omega))^m}
\leqslant Ce^{K/T^p}\|y^0\|_{L^2(\Omega)},
\end{equation*}
where $K>0$  is the constant given by Proposition \ref{prop ine obs 2} with $N=2R$ and $p\geqslant 2$ is given in Proposition \ref{prop ine obs 1}.

Letting the system evolve freely a little bit if needed, we may assume without loss of generality that $y^0-\overline y^0\in H^1_0(\Omega)$. Indeed, by the regularizing effect, it is very easy to deduce that  for any solution $(\overline{y},\overline{u})$ to \eqref{syst linearise}, 
there exists some $C(T)>0$ such that for any solution $(y,0)$  to  \eqref{syst primal} on $[0,\frac{T}{2}]$, we have $y\left(\frac{T}{2}\right)-\overline y\left(\frac{T}{2}\right)\in H^1_0(\Omega)$ and   $$||y\left(\frac{T}{2}\right)-\overline y\left(\frac{T}{2}\right)||_{H^1(\Omega)}\leqslant C(T)||y^0-\overline y^0||_{L^2(\Omega)}.$$
%Moreover, $C(T,||u||_{L^\infty((0,T),W^{1,\infty}(\Omega))})$ can be chosen to be continuous with respect to $||u||_{L^\infty((0,T),W^{1,\infty}(\Omega))}$.
Hence, if $||y^0-\overline y^0||_{L^2(\Omega)}$ is small, so is $||y\left(\frac{T}{2}\right)-\overline y\left(\frac{T}{2}\right)||_{H^1(\Omega)}$, so that the condition \eqref{cieta} is sufficient for our argument to be valid.

Following \cite[p. 24]{MR3023058}, we introduce the cost of controllability given by
$$\gamma(t)= Ce^{K/t^p},\,\,\,t\in (0,T),$$ and the following weight functions 
$$\rho_{F}(t)=e^{-\frac{\alpha}{(T-t)^{p+1}}}, t\in [0,T]$$
and 

$$\rho_0(t)=e^{\frac{K}{((q-1)(T-t))^{p}}-\frac{\alpha}{q^{2p+2}(T-t)^{p+1}}},\,t\in\left[T\left(1-\frac{1}{q^2}\right),T\right],$$
extended on $[0, T\left(1-\frac{1}{q^2}\right)]$ by 
$$\rho_0(t)=\rho_0\left(T\left(1-\frac{1}{q^2}\right)\right), \,t\in [0, T\left(1-\frac{1}{q^2}\right)],$$
for some parameters  $q>1$ and $\alpha>0$ to be chosen later on.

We remark that $\rho_F$ and $\rho_0$ are non-increasing, that they verify $\rho_F(T)=\rho_0(T)=0$, and are related by 
$$\rho_0(t)=\rho_F(q^2(T-t)+T)\gamma((q-1)(T-t)),\,\, t\in\left[T\left(1-\frac{1}{q^2}\right),T\right].$$ 

We introduce for some $\beta>0$ the weight function
$$\rho(t)=e^{-\frac{\beta}{(T-t)^{p+1}}}.$$
We remark that
\begin{equation*}%\label{cpoi}
\rho_F\leqslant C\rho,\,\,\mbox{} \rho_0\leqslant C\rho, \,\, \mbox{} |\rho'|\rho_0\leqslant C\rho^2,
\end{equation*}
as soon as $\beta>0$ is chosen small enough, precisely 
\begin{equation}\label{condbeta} \beta<\frac{\alpha}{q^{2p+2}}.
\end{equation}
We introduce the following spaces:
$$
\mathcal F=\{f\in  L^2((0,T)\times\Omega), \frac{f}{\rho_F}\in   L^2((0,T)\times\Omega)\},
$$
%$$\mathcal U=\{u\in L^\infty((0,T), W^{1,\infty}(\Omega))^m \\ \mbox{ such that }\frac{u}{\rho_0}\in L^\infty((0,T),W^{1,\infty}(\Omega))^m\},$$

\textcolor{black}{$$\mathcal U=\{u\in L^2((0, T ) \times\Omega),\frac{u}{\rho_0}\in  L^{\infty}((0, T), W^{1,\infty})\}$$
and
$$\mathcal Z=\{z\in L^2((0, T ) \times\Omega),\frac{z}{\rho}\in H^1((0, T), L^2)\cap L^2((0,T),H^2\cap H^1_0)\},$$}
%\begin{multline*}\mathcal Z=\{z\in C^0([0,T],H^1_0(\Omega))\cap L^2((0,T),H^{2}(\Omega)\cap H^1_0(\Omega))\cap H^{1}((0,T),L^{2}(\Omega))\\ \mbox{ such that }\frac{z}{\rho}\in   C^0([0,T],H^1_0(\Omega))\cap L^2((0,T),H^{2}(\Omega)\cap H^1_0(\Omega))\cap H^{1}((0,T),L^{2}(\Omega))\}, 
%\end{multline*}
endowed with the weighted Sobolev norms naturally induced by the definition of these spaces.

Following  \cite[Proofs of Propositions 2.5, 2.8]{MR3023058}  in the spirit of   \cite[Section 7.2 and Appendix 5]{klb}, it is easy to obtain the following result.
\begin{prop}\label{pfin}
For any $z^0\in H^1_0(\Omega)$  and any $f\in \mathcal F$,
 there exists $v\in \mathcal U$ such that the solution $z$ of 
\begin{equation*}%\label{syst lin primal}
\left\{\begin{array}{rll}
\partial_tz&=\Delta z+\Div(\overline u z)+\Div(\theta\overline y Bv)+f&\mbox{in } Q_T,\\
z&=0&\mbox{on } \Sigma_T,\\
z(0,\cdot)&=z^0&\mbox{in }\Omega,
\end{array}\right.
\end{equation*}  
verifies $z\in\mathcal Z$ (and hence $z(T)=0$). 
%\textcolor{red}{pourquoi ne pas mettre l'estimation?}
%Moreover, one has 
%\begin{equation}\label{complique}\left\|\frac{z}{\rho}\right\|_{C^0([0,T],H^1_0(\Omega))\cap L^2((0,T),H^{2}(\Omega))\cap H^{1}((0,T),L^{2}(\Omega))}\leqslant C\left(||z^0||_{H^{1}(\Omega)}+ \left\|\frac{f}{\rho_F}\right\|_{L^2((0,T)\times\Omega)}\right).\end{equation}

%verifies 
%\begin{equation}\label{okpourci}\left\|\frac{z}{\rho}\right\|_{ L^2((0,T),H^{2}(\Omega))\cap H^{1}((0,T),L^{2}(\Omega))}\leqslant C\left(||z^0||_{H^{1}(\Omega)}+ \left\|\frac{f}{\rho_F}\right\|_{L^2((0,T)\times\Omega)}\right).\end{equation}

\end{prop}
%\textbf{Proof of Proposition \ref{pfin}.}
%\textcolor{black}{a faire}
%\cqfd

To conclude, we use the following inverse mapping theorem:

\begin{theo}[see \cite{alekseev_fonc_impl87}]\label{fonc impl}
Let $\mathcal{X}$ and $\mathcal{Y}$ be Banach spaces and let $M :  \mathcal{X}\mapsto \mathcal{Y}$ be a $\mathcal{C}^1$ mapping. 
\textcolor{black}{Consider $x_0\in \mathcal{X}$ and $y_0:=M(x_0)\in \mathcal{Y}$.}
Assume that the derivative
$\textcolor{black}{M' (x_0 )}%M'(0)
: \mathcal{X} \mapsto \mathcal{Y}$ is onto. 
Then, there exist $\eta > 0$, a mapping 
$W : B_{\eta}(y_0)\subset \mathcal{Y}\mapsto \mathcal{X}$ and a constant $K > 0$
satisfying:
\begin{equation*}
 \left\{
\begin{array}{l}
 W(z)\in \mathcal{X}\mbox{ and } M(W(z))=z~\forall z\in B_{\eta}(y_0),\\
 \|W(z)-x_0\|_{\mathcal{X}}\leqslant K\|z-y_0\|_{\mathcal{Y}}~\forall z\in B_{\eta}(y_0).
\end{array}
\right.
 \end{equation*}
\end{theo}

\textbf{Proof of Theorem \ref{theo: contr traj 2}.}
%Our proof is close from \cite[Section 8]{klb}. 
We are looking for a solution in the form
\begin{equation*}
%\label{chang trajec}
 y(x,t)=\overline{y}(x,t)+w(x,t),~u(x,t)=\overline{u}(x,t)+\theta(x)Br(x,t),
\end{equation*}
where $(y,u)$ and $(\overline{y},\overline{u})$ are solution to the Systems \eqref{syst primal} and \eqref{syst linearise}, respectively.
Then $(w,r)$ is solution to 
\begin{equation*}
 \left\{\begin{array}{rcll}
         N(w,r)&:=&\partial_tw-\Delta w-\Div(\overline{u}w+\theta Br\overline{y}+\theta Brw)=0
         &\mbox{in } Q_T,\\
       w&=&0&\mbox{on } \Sigma_T,\\
	w(0,\cdot)&=&y^0-\overline y^0 &\mbox{in } \Omega.\\
	%w(T,\cdot)&=&0&\mbox{in } \Omega.
        \end{array}
\right.
\end{equation*}

We introduce the following spaces:
$${\mathcal{X}}:=\{(w,r)\in \mathcal Z\times\mathcal U\mbox{ such that } \partial_tw-\Delta w-\Div(\overline{u}w+\theta Br\overline{y})\in \mathcal F\},$$
endowed with the norm 
$$||(w,r)||_{{\mathcal X}}=||w||_{\mathcal Z}+||r||_{\mathcal U}+ ||\partial_tw-\Delta w-\Div(\overline{u}w+\theta Br\overline{y})||_{\mathcal F},$$ and the space 
$${\mathcal Y}=\mathcal{F}\times H^1_0(\Omega),$$
endowed with the norm 
$$||(f,z^0)||_{\mathcal Y}:=||f||_{\mathcal F}+||z^0||_{H^1(\Omega)}.$$
Introduce the mapping $M$ given by
\begin{equation*}
\begin{array}{cccc}
 M:&{\mathcal{X}}&\rightarrow& \mathcal{Y}\\
 &(w,r)&\mapsto&(N(w,r),w(0,\cdot)).
\end{array}
\end{equation*}
Let us determine what are the conditions on $q,\alpha,\beta$ ensuring that $M$ is well-defined. It is clear that 
$$||w(0,\cdot)||_{H^1(\Omega)}\leqslant ||w||_{C^0([0,T],H^1_0(\Omega))}\leqslant C\left\|\frac{w}{\rho}\right\|_{C^0([0,T],H^1_0(\Omega))}  \leqslant ||(w,r)||_{\mathcal X}.$$

Now, we remark that by definition of the space $\mathcal X$, we have 

$$ ||\partial_tw-\Delta w-\Div(\overline{u}w+\theta Br\overline{y})||_{\mathcal F}\leqslant   ||(w,r)||_{\mathcal X}.$$

Hence, the only difficulty is to treat the bilinear part $\Div(\theta wBr)$. 
We remark that 
$$\displaystyle\left\|\frac{\Div(\theta wBr)}{\rho_F}\right\|_{L^2((0,T)\times\Omega)}\leqslant C\displaystyle\left\|\frac{r}{\rho_F^{\frac{1}2}}\right\|_{L^\infty((0,T),W^{1,\infty}(\Omega))}\left\|\frac{w}{\rho_F^{\frac{1}2}}\right\|_{L^2((0,T),H^1(\Omega))}.$$
We can impose that $\rho^2\leqslant C \rho_F$ and $\rho_0^2\leqslant C \rho_F$ as soon as 
\begin{equation}\label{condbeta2}\alpha<2\beta \,\,\mbox{ and } q^{2p+2}< 2.
\end{equation}
Remark that these conditions are compatible with condition \eqref{condbeta}.

Hence, under conditions \eqref{condbeta} and \eqref{condbeta2}, we deduce that 

$$\begin{aligned}\left\|\frac{\Div(\textcolor{black}{\theta  wBr})}{\rho_F}\right\|_{ L^2((0,T)\times\Omega)}&\leqslant C\displaystyle\left\|\frac{r}{\rho_0}\right\|_{L^\infty((0,T),W^{1,\infty}(\Omega))}\left\|\frac{w}{\rho}\right\|_{L^2((0,T),H^1(\Omega))}\\&\leqslant C||(w,r)||_{\mathcal X}^2.\end{aligned}$$

We conclude that under these conditions, $M$ is indeed well-defined and continuous. Moreover,
we remark that $M(0,0)=(0,0)$ and $M$ is of class $\mathcal{C}^1$ as a sum of a  %linear continuous 
\textcolor{black}{continuous linear} function and a %quadratic continuous 
\textcolor{black}{continuous quadratic} function.
Furthermore, Proposition \ref{pfin} exactly means that 
$M'(0,0)$ is onto (see Remark \ref{eqt}), \textcolor{black}{when 
\begin{equation}\label{eq:cond eta}
\frac{\alpha}{q^{2p+p}}<1
\end{equation}
and $\eta\in(0,1)$ is chosen as 
$$\eta:=\frac{\alpha}{q^{2p+p}}.$$
Conditions \eqref{condbeta}, \eqref{condbeta2} and \eqref{eq:cond eta} can be summarized as follows:
$$\frac{\alpha}{2}<\beta<\frac{\alpha}{q^{2p+p}}<1 \mbox{ and }q^{2p+2}<2,$$
which is satisfied for $q=(3/2)^{1/(2p+2)}> 1$ (remind that $p\geqslant 2$, so that $2p+p\geqslant 2p+2\geqslant 1$), $\alpha=1$ and $\beta=7/12$.}
 Theorem \ref{fonc impl} leads to the conclusion.

\cqfd

\section{Example of a non-controllable trajectory with a reduced number of controls}\label{sec:CE}

In this section, we give %an example of trajectory 
\textcolor{black}{the example of a trajectory} which does not satisfy condition \eqref{cdet}  and for which the local controllability to the trajectories does not hold.

Consider $\overline{u}\in L^{\infty}(Q_T)^m$ \textcolor{black}{ which is independent of the time variable} and will be determined later on.
Assume that for each $\overline{y}^0\in L^2(\Omega)\setminus \{0\}$ the following system is locally controllable to the trajectories with a control operator $B$ \textcolor{black}{to be chosen later on:}
\begin{equation}\label{syst linearisece}
\left\{\begin{array}{lll}
\partial_t\overline{y}&=\Delta \overline{y}+\Div(\overline{u} \overline{y})&\mbox{in } Q_T,\\
\overline{y}&=0&\mbox{on } \Sigma_T,\\
\overline{y}(0,\cdot)&=\overline{y}^0&\mbox{in }\Omega.
\end{array}\right.
\end{equation}
Then, for any $\varepsilon\in(0,1)$, there exists $u\in L^{\infty}(Q_T)^m$ such that
\begin{equation*}%\label{syst primal}
\left\{\begin{array}{lll}
\partial_ty&=\Delta y+\Div(u y)&\mbox{in } Q_T,\\
y&=0&\mbox{on }\Sigma_T,\\
y(0,\cdot)&=(1-\varepsilon)\overline{y}^0&\mbox{in }\Omega,\\
y(T,\cdot)&=\overline{y}(T)&\mbox{in }\Omega,
\end{array}\right.
\end{equation*}
where $u=\overline{u}+Bv$ with $\Supp(v)\subset (0,T)\times\omega$.
We remark that $(z,w):=(y-\overline{y},yv)$ is solution to
\begin{equation}\label{syst z}
\left\{\begin{array}{lll}
\partial_tz&=\Delta z+\Div(\overline{u} z)+\Div(Bw)&\mbox{in } Q_T,\\
z&=0&\mbox{on }\Sigma_T,\\
z(0,\cdot)&=\varepsilon\overline{y}^0&\mbox{in }\Omega,\\
z(T,\cdot)&=0&\mbox{in }\Omega.
\end{array}\right.
\end{equation}
We deduce that \textcolor{black}{the linear control} system \eqref{syst z} is null controllable at time $T>0$, then approximately controllable at time $T>0$.
It is well known that the approximate controllability of System \eqref{syst z} on  $(0,T)$  implies the following property, called the Fattorini-Hautus test
(see \textit{e.g.} \cite{olive_bound_appr_2014}) : 
%(see \cite[Theorem 1 \& Section 3]{olive_bound_appr_2014}):
%\begin{theo}\label{theo:fattorini2}
%System \eqref{syst primal} is approximately controllable on  $(0,T)$, if and only if 
for every $s\in\mathbb{C}$ and every $\varphi \in H^2(\Omega)\cap H^1_0(\Omega)$, 
\begin{equation}\label{prop fat}
\left. \begin{array}{ll}
-\Delta  \varphi-\overline{u}\cdot\nabla\varphi  &= s\varphi \mbox{ in } \Omega\\\noalign{\smallskip}
B^*\nabla \varphi&= 0\mbox{  in }\omega
\end{array}\right\}
\Rightarrow \varphi = 0.
\end{equation}

Now, we give an explicit situation  in contradiction with \eqref{prop fat}.
Consider $\Omega=(0,\pi)^2$, $B^*=(1,0)$ and  \textcolor{black}{$s=25$} ($\omega$ and $\overline u$ will be chosen later on).
The goal is to find a nontrivial solution $\varphi \in H^2(\Omega)\cap H^1_0(\Omega)$ of

\begin{equation}\label{prop fat 2}
\left\{ \begin{array}{ll}
-\Delta  \varphi -\overline u \cdot \nabla \varphi&\textcolor{black}{=25} \varphi \mbox{ in } \Omega,\\\noalign{\smallskip}
\partial_{x_1} \varphi&= 0\mbox{  in }\omega.
\end{array}\right.
\end{equation}

We introduce two functions $f$ and $g$ defined on $\Omega$ and  given by 
$$\textcolor{black}{f(x_1,x_2)=\sin(3x_1)\sin(4x_2),\,g(x_1,x_2)=\frac{2\sqrt{2}}{5}\left (-\sin(5x_2)+\cos(5x_2)\right).}$$
Remark that $f\in  C^\infty(\Omega)$ and $g\in C^\infty(\Omega)$ are chosen in such a way that 
\begin{equation}
\label{fgpro}
\left\{ \begin{array}{ll}
-\Delta f &= \textcolor{black}{25}f \mbox{ on } \Omega,\\
-\Delta g &= \textcolor{black}{25}g \mbox{ on } \Omega,\\
\partial_{x_1}g&=0 \mbox{ on } \Omega,\\
f&=0 \mbox{ on } \partial\Omega,\\
f\left(\frac{\pi}{2},\frac{\pi}{4} \right )&=g\left(\frac{\pi}{2},\frac{\pi}{4} \right ),\\
\nabla f\left(\frac{\pi}{2},\frac{\pi}{4} \right )&=\nabla g\left(\frac{\pi}{2},\frac{\pi}{4} \right ).\\
\end{array}\right.
\end{equation}

Now, let us consider some cut-off function $\chi\in C^\infty_0(\mathbb R^2)$ such that $\chi\in[0,1]$,\, $\chi \equiv 0$ on $\mathbb R\setminus[\frac{\pi}{4},\frac{3\pi}{4}]\times \mathbb R\setminus[\frac{\pi}{8},\frac{3\pi}{8}]$ and $\chi=1$ on $[\frac{3\pi}{8},\frac{5\pi}{8}]\times [\frac{3\pi}{16},\frac{5\pi}{16}]$.
 For a parameter $h\in (0,1)$, we call 
$$\chi_h(x_1,x_2):=\chi\left(\frac{\pi}{2}+\frac{x_1-\pi/2}{h},\frac{\pi}{4}+\frac{x_2-\pi/4}{h} \right).$$
Note that $\chi_h$ is supported in 
$$V_h:=\left[\frac{\pi}{2}-\frac{h\pi}{4},\frac{\pi}{2}+\frac{h\pi}{4} \right ]\times \left[\frac{\pi}{4}-\frac{h\pi}{8},\frac{\pi}{4}+\frac{h\pi}{8} \right ],$$
verifies
\begin{equation}\label{ouv1}
\chi_h=1\,\, \mbox{ on }\, \left[\frac{\pi}{2}-\frac{h\pi}{8},\frac{\pi}{2}+\frac{h\pi}{8} \right ]\times \left[\frac{\pi}{4}-\frac{h\pi}{16},\frac{\pi}{4}+\frac{h\pi}{16} \right ]:=W_h
\end{equation}
and 
\begin{equation}\label{ech}|\chi_h|\leqslant 1,\,|\partial_{x_2} \chi_h|\leqslant \frac{C}{h}\mbox{ on } \mathbb R^2,
\end{equation}
for some $C>0$ independent on $h$.
Now, we introduce 
\begin{equation}\label{dvp}\varphi_h=\chi_h g+(1-\chi_h)f.\end{equation}
We remark that for any $h\in (0,1)$, $\varphi_h \not \equiv 0$ since it coincides with $f$ outside $V_h\not = \emptyset$ and with $g$ on $W_h\not=\emptyset$.
Moreover, one has
\begin{equation}\label{vphi}\partial_{x_2}\varphi_h= \partial_{x_2} \chi_h(g-f)+\chi_h\partial_{x_2} g + (1-\chi_h)\partial_{x_2}f.\end{equation}
By the two last lines of \eqref{fgpro} and  Taylor expansions, for $h\in(0,1)$, we have 
\begin{equation}\label{fpg} |f-g|\leqslant Ch^2 ,\, |\partial_{x_2}f-\partial_{x_2}g|\leqslant Ch \mbox{ on } V_h,
\end{equation}
for some $C>0$ independent on $h$.
From \eqref{ech}, \eqref{vphi} and \eqref{fpg}, we deduce that 
$$|\partial_{x_2} \varphi_h-\partial_{x_2}g|\leqslant   |\partial_{x_2} \chi_h|\,|g-f|+|1-\chi_h|\,|\partial_{x_2} g-\partial_{x_2}f|\leqslant Ch\mbox{ on } V_h,$$
for some $C>0$ independent on $h$.
Since $$\partial_{x_2}g\left(\frac{\pi}{2},\frac{\pi}{4} \right)=\textcolor{black}{4}>0,$$
we deduce that there exists $h_0>0$ small enough such that, $\partial_{x_2} \varphi_{h_0} \geqslant C$ on $V_{h_0}$ for some $C>0$.
Accordingly to \eqref{ouv1}, we choose 
$\omega=W_{h_0},$
and 
\begin{equation*}
\overline{u}:=\left\{\begin{array}{ll}(0,0)&\mbox{ in } (0,\pi)^2\setminus (V_{h_0}\setminus W_{h_0}),
\\(0,-\frac{\textcolor{black}{25}\varphi_{h_0}+\Delta \varphi_{h_0}}{\partial_{x_2}\varphi_{h_0}})&\mbox{ otherwise}.
 \end{array}\right .
\end{equation*}

Remark that $\overline u$ is well defined: by construction, $\partial_{x_2}\varphi_{h_0}\geqslant C$ where $-\Delta \varphi_{h_0}-\textcolor{black}{25}\varphi_{h_0} \not =0$ (which is included in $V_{h_0}\setminus W_{h_0}$).

Moreover, $\overline u$ is of class $C^\infty$ on $\Omega$.
To conclude, we remark that by \eqref{fgpro} and \eqref{dvp},
\begin{equation*}
\left\{\begin{array}{ll}-\Delta  \varphi_{h_0}-\overline{u}\cdot\nabla\varphi_{h_0}  &=\textcolor{black}{ 25}\varphi_{h_0}\mbox{ in }\Omega,\\
\varphi_{h_0}&=0\mbox{ in }\partial\Omega,\\
\partial_{x_1}\varphi_{h_0}&=0\mbox{ in }\omega,\\
\varphi_{h_0}&\neq0.\end{array}\right.
\end{equation*}

Hence, we obtain a contradiction with the Fattorini-Hautus test, which concludes our proof.
\color{black}

\vspace{0.2cm}

\textbf{{\Large Acknowledgements}}

\vspace*{0.1cm}

\textcolor{black}{The authors would like to thank the anonymous referee for his/her comments that allowed us to correct some errors and imprecisions and to simplify the exposition of some proofs and assumptions, thus improving the overall presentation of the paper.}

\footnotesize
\bibliographystyle{plain}
\bibliography{biblio.bib}

\end{document}